\documentclass[12pt]{article}
\usepackage{latexsym}
\usepackage{amsfonts}
\usepackage{amsmath}
\usepackage{latexsym}
\usepackage{amsfonts}
\usepackage{amssymb}
\usepackage{rawfonts}
\input{prepictex}
\input{pictex}
\input{postpictex}

\def\bcp{\mathbb C \mathbb P}
\def\eea{\end{eqnarray*}}
\def\hull{\operatorname{\bf Hull} }

\newtheorem{main}{Theorem}

\newtheorem{thm}{Theorem}[section]
\newtheorem{lem}[thm]{Lemma}
\newtheorem{prop}[thm]{Proposition}

\newtheorem{defn}[thm]{Definition}

\newenvironment{proof}{\medskip \noindent
{\bf Proof.}}{\hfill \rule{.5em}{1em}
\\}
\newenvironment{proofA}{\medskip \noindent
{\bf Proof of Theorem A.}}{\hfill \rule{.5em}{1em}
\\}

\newenvironment{xpl}{\mbox{ }\\ {\bf  Example}\mbox{ }}{
\hfill $\diamondsuit$\mbox{}\bigskip}

\def\ZZ{{\mathbb Z}}
\def\RR{{\mathbb R}}

\def\span{\mbox{\bf span }}
\begin{document}
\sloppy
\title{Four-Manifolds, Curvature Bounds,\\ and Convex Geometry}

\author{Claude LeBrun\thanks{Supported 
in part by  NSF grant DMS-0604735} 
\\ 
SUNY Stony
 Brook 
  }

\date{November 14, 2006}

\maketitle

\section{Introduction} 
Seiberg-Witten theory leads to a remarkable family of curvature estimates governing
the Riemannian geometry of compact  $4$-manifolds, and these, for example, 
imply interesting results concerning 
the existence and/or uniqueness 
of Einstein metrics on such spaces. 
The primary purpose of the present article is to introduce a simplified, user-friendly  
repackaging of the information conveyed by the Seiberg-Witten equations
into a single, easily understood numerical invariant  that appears to play the starring r\^ole 
in the relevant curvature estimates. In addition, this article
contains some new results concerning boundary cases of the curvature estimates that 
strengthen what was previously known. 

The gist of the matter can
 be summarized as follows.  Suppose that $M$ is a smooth 
 compact oriented $4$-manifold with $b_+(M) \geq 2$. By considering a geometrically rich 
system of PDE called the 
Seiberg-Witten equations, one may then define a certain finite subset ${\mathfrak C}\subset 
H^2(M, \RR)$  that  depends only on the orientation and smooth structure 
of $M$. The elements of ${\mathfrak C}$ are called {\em monopole classes},
and   are, by definition,  the first Chern classes  of those spin$^c$ structures
on $M$ 
for which the  the Seiberg-Witten equations have solutions for all metrics.
Now, while   the elements of ${\mathfrak C}$ are all {\em integer} classes, 
we wish to focus here on the fact that ${\mathfrak C}$ sits in  a 
real vector space, as this allows us to consider its 
 convex hull  $\hull({\mathfrak C})$. 
Because  ${\mathfrak C}$ is finite, $\hull({\mathfrak C})$ is 
necessarily compact. We can therefore define a real-valued 
invariant of $M$ by setting 
$$
\beta^2 (M) = \max  ~\left\{~ {\mathbf v}^2 ~|~ {\mathbf v}\in \hull({\mathfrak C})~\right\}
$$
when ${\mathfrak C}\neq \varnothing$, while  setting $\beta^2(M)=0$ if  
${\mathfrak C}= \varnothing$. Here 
${\mathbf v}^2=\langle {\mathbf v}\smallsmile {\mathbf v} , [M]\rangle$ 
denotes the {\em intersection pairing} of a class 
${\mathbf v}\in H^2(M, \RR )$ with itself.  Because $0\in \hull ({\mathfrak C})$ whenever 
${\mathfrak C}\neq \varnothing$, one automatically has 
$\beta^2 (M) \geq 0$; but, more importantly, there are actually 
many $4$-manifolds $M$ for which  $\beta^2 (M) > 0$.
It is this last fact  that gives the following result most of its  interest: 

\begin{main} \label{summa} 
Let $M$ be a compact oriented $4$-manifold with $b_+\geq 2$. 
Then any metric $g$ on $M$ satisfies the curvature estimates 
\begin{eqnarray}
\int_M s^2d\mu &\geq& 32\pi^2 \beta^2(M)\label{mocha}\\
\int_M (s-\sqrt{6}|W_+|)^2d\mu &\geq& 72\pi^2 \beta^2(M)\label{java}
\end{eqnarray}
where $s$ and $W_+$ respectively denote the scalar and Weyl curvatures of $g$. 
Moreover, if $M$ carries a non-zero monopole class, 
equality  occurs in either (\ref{mocha}) or (\ref{java})   
if and only if  $g$ is K\"ahler-Einstein, with negative Einstein constant. 
\end{main}

Now,  in an important  respect,  Theorem \ref{summa}  is ostensibly  weaker than a 
result that the author has published elsewhere \cite{lebsurv2}. Indeed, as we shall 
see below, there is  a `softer' invariant, called  $\alpha^2(M)$,  that can 
be defined in terms of  ${\mathfrak C}(M)$ via a   complicated
minimax process; and na\"{\i}ve comparison of the definitions of $\alpha^2$ and $\beta^2$ would  
lead one  merely to expect that 
$$\beta^2(M) \leq \alpha^2 (M) \leq b_+(M) \beta^2(M).$$
Yet  \cite{lebsurv2}  inequalities  (\ref{mocha}) and  (\ref{java}) can still be shown to hold  
even when $\beta^2 (M)$
is replaced by $\alpha^2 (M)$, yielding  an apparently stronger statement. 
Oddly enough, however, it seems  that  {\em in practice}   one consistently 
has $$\alpha^2(M)= \beta^2 (M),$$ 
so that  modifying  (\ref{mocha}) or (\ref{java}) in this manner 
effectively  seems to 
yield no added punch. 
The fact that $\alpha^2$ and $\beta^2$  typically  coincide will only  partially be explained
here, 
via some simple results of distinctly limited scope.
But the upshot is that the finite configuration   ${\mathfrak C}\subset 
H^2(M, \RR)$ consistently displays  some unanticipated  geometrical properties 
that ought to be understood  more precisely. 

In yet a  different direction,  Theorem \ref{summa}  contains some 
essentially new geometric information, because  the stated  characterization 
of the equality case of 
  (\ref{java}) was not previously known. 
 The issue boils down to a   problem in almost-K\"ahler 
  geometry, and is   eventually resolved by  Theorem \ref{newnew}.
  
 Finally, it should be pointed out that  
  the convex hull of the set of monopole
 classes first   appeared   in the  context of
 $3$-manifold theory, where Kronheimer and Mrowka \cite{KM3} 
 used it to provide a  new characterization of the Thurston norm.
 Although these earlier results ultimately have  little bearing on 
 the ideas developed here, they undoubtedly exerted a 
  powerful subconscious 
 influence   on  the conceptualization of the present work.
 The author would therefore like to express his indebtedness   to Kronheimer and Mrowka
by drawing the reader's attention to their deep and beautiful   paper.

\section{Rudiments of $4$-Dimensional Geometry}
  
This article will make frequent  reference to  a constellation of basic facts 
regarding  $4$-dimensional   geometry   which, though largely   familiar to the  
 {\em cognoscenti},  would completely confuse  the neophyte if left unexplained.
  For clarity's sake, we will  therefore begin with a  quick review of  the key points.

Many peculiar features of  $4$-dimensional geometry are directly attributable  the fact that 
the bundle  of 2-forms over an
oriented Riemannian $4$-manifold $(M,g)$  invariantly  decomposes 
as the direct sum
\begin{equation} 
\Lambda^2 = \Lambda^+ \oplus \Lambda^- 
\label{deco} 
\end{equation}
 of the eigenspaces of 
the Hodge star 
operator
$$\star: \Lambda^2 \to \Lambda^2.$$
The sections of $\Lambda^+$ are   called {\em self-dual 
2-forms}, while the  
sections of $\Lambda^-$   are  called {\em anti-self-dual
2-forms}. The decomposition (\ref{deco}) is, moreover,  {\em conformally invariant},
in the sense that it is left unchanged if $g$ is multiplied by an arbitrary smooth positive function. 
An arbitrary $2$-form can thus be uniquely expressed as 
$$\varphi = \varphi^+ + \varphi^-,$$
where $\varphi^\pm \in \Lambda^\pm$, and we then have
$$\varphi \wedge \varphi = \Big( |\varphi^+|^2 - |\varphi^-|^2\Big) d\mu_g , $$
where $d\mu_g$ denotes the metric volume form associated with the
fixed orientation. 

The decomposition (\ref{deco})  in turn leads to a decomposition of the 
 Riemann curvature tensor. Indeed, identifying the curvature tensor
 of $g$ with the  self-adjoint linear map 
 \begin{eqnarray*}
{\mathcal R} : \Lambda^2 &\longrightarrow& \Lambda^2\\
\varphi_{jk}&\longmapsto&{\textstyle \frac{1}{2}}\varphi_{\ell m} {R^{\ell m}}_{jk}
\end{eqnarray*}
we obtain a decomposition
 \begin{equation}
\label{curv}
{\mathcal R}=
\left(
\begin{array}{c|c}&\\
W_++\frac{s}{12}&\mathring{r}\\ 
&\\
\hline 
&\\
\mathring{r} & W_-+\frac{s}{12}\\
&\\
\end{array} \right) . 
\end{equation}
where $W_++\frac{s}{12} : \Lambda^+\to \Lambda^+$, etc. 
Here $W_+$ is the trace-free piece of its block, and is the called
the {\em self-dual Weyl curvature} of $(M,g)$; the anti-self-dual
Weyl curvature $W_-$ is defined analogously.
Both of the objects are conformally invariant, in the sense that the tensors 
${(W_\pm)^j}_{k\ell m}$ both remain unaltered if $g$ is multiplied by any
smooth positive function. Note that 
the  {scalar curvature} $s$ is understood to act in (\ref{curv}) by scalar multiplication,
while  the  trace-free  Ricci curvature
$\mathring{r}_{jk}={R^\ell}_{j\ell k}-\frac{s}{4}g_{jk}$ 
acts on 2-forms by
$$\varphi_{jk} \mapsto ~ 
 \mathring{r}_{\ell [j}{\varphi^{\ell}}_{k]}.$$
  
Next,  let us  suppose that  $(M,g)$ is a {\em  compact} oriented
Riemannian  $4$-manifold.
The  Hodge theorem then tells us that every de Rham class
on $M$ has a unique harmonic representative. In particular, 
 we therefore have a canonical identification 
$$H^2(M,{\mathbb R}) =\{ \varphi \in \Gamma (\Lambda^2) ~|~
d\varphi = 0, ~ d\star \varphi =0 \} .$$
However,  the Hodge star operator $\star$ defines an involution of 
the right-hand side, giving rise to  an eigenspace decomposition
\begin{equation}
	H^2(M, {\mathbb R}) = {\mathcal H}^+_{g}\oplus {\mathcal H}^-_{g},
	\label{harm}
\end{equation}
where
$${\mathcal H}^\pm_{g}= \{ \varphi \in \Gamma (\Lambda^\pm) ~|~
d\varphi = 0\} $$
are the spaces of self-dual and anti-self-dual harmonic forms.
The intersection form 
\begin{eqnarray*}
H^{2}(M, {\mathbb R})\times H^{2}(M, {\mathbb R})	
 & \longrightarrow & ~~~~ {\mathbb R}  \\
	( ~ {\mathbf  b}~ , ~{\mathbf c} ~) ~~~~~~~~~ & 
	\longmapsto  & {\mathbf b} \cdot  {\mathbf c} := 
	\langle {\mathbf b} \smallsmile {\mathbf c}, [M]\rangle
\end{eqnarray*}
becomes positive-definite when restricted to 
${\mathcal H}^+_{g}$, and  negative-definite when restricted to 
${\mathcal H}^-_{g}$. Moreover, these two subspaces are 
mutually orthogonal with respect to the intersection form. 
The numbers
$b_{\pm} (M)= \dim {\mathcal H}^\pm_{g}$ are  therefore oriented 
homotopy invariants of $M$. 
Their difference 
$$\tau (M) = b_+(M) -b_-(M)$$
is called the {\em signature} of $M$. By the Hirzebruch Signature Theorem,
it coincides with $\langle\frac{1}{3}p_1(TM) , [M]\rangle$, and so can be 
expressed as a curvature integral 
 \begin{equation}
\tau (M)= \frac{1}{12\pi^2}\int_M \Big(|W^+|^2
- |W^-|^2  \Big)d\mu 
\end{equation}
for  any Riemannian metric $g$ on $M$.
This, of course,  is  analogous to the generalized Gauss-Bonnet formula 
 \begin{equation}
\chi (M)= \frac{1}{8\pi^2}\int_M \left(\frac{s^2}{24}+|W^+|^2
+ |W^-|^2  -\frac{|\stackrel{\circ}{r}|^2}{2}
\right)d\mu 
\end{equation}
for the Euler characteristic. 

\begin{lem} \label{harmless} 
Let $\psi$ be a closed $2$-form on a compact oriented Riemannian $4$-manifold
$(M,g)$.
Let ${\mathbf v}= [\psi ]$ denote the de Rham class of $\psi$, and 
use \eqref{harm} to write 
$${\mathbf v} = {\mathbf v}^+ + {\mathbf v}^-$$
where ${\mathbf v}^\pm \in {\mathcal H}^\pm_g$.
Then 
$$\int_M |\psi^+|^2 d\mu_g \geq ({\mathbf v}^+)^2,$$
with equality iff $\psi^+$ is a harmonic form.
\end{lem}
\begin{proof}
Let $\phi$ be the unique harmonic form cohomologous to $\psi$. 
Since $\phi$ is then the de Rham representative of ${\mathbf v}$ of minimal  $L^2$-norm,
we therefore have 
$$
\int_M (|\psi^+|^2+|\psi^-|^2) d\mu \geq \int_M (|\phi^+|^2+|\phi^-|^2) d\mu ~,
$$
with equality iff $\psi = \phi$. However, 
$$
\int_M (|\psi^+|^2-|\psi^-|^2) d\mu = \int_M (|\phi^+|^2-|\phi^-|^2) d\mu ~,
$$
since  $\int \psi \wedge \psi= \int \phi\wedge \phi$ by Stokes' theorem. 
Averaging these expressions, we therefore have
$$
\int_M |\psi^+|^2d\mu \geq \int_M |\phi^+|^2 d\mu = \int_M \phi^+\wedge \phi^+
 = ({\mathbf v}^+)^2~,
$$
with equality iff 
$ \psi^+= (\psi + \star\psi)/2$ is closed. 
\end{proof}

When using this result, it is important to  remember that the 
decomposition
\eqref{harm} depends on the metric $g$, as consequently
does  the number $({\mathbf v}^+)^2$. This  makes it
vital  to better understand 
the natural map 
\begin{eqnarray*}
\{\mbox{metrics on }M\}  &\longrightarrow& Gr^+_{b_+}[H^2(M,\RR)]\\
g~~~~&\longmapsto& ~~~~{\mathcal H}_g^+
\end{eqnarray*}
from the infinite-dimensional space of all metrics to 
the finite-dimensional Grassmannian of $b_+(M)$-dimensional 
subspaces of $H^2(M,\RR)$ on which the intersection form is 
positive-definite. This map  is called the
{\em period map} of $M$. It is easily seen to be  invariant under 
both conformal rescaling and 
the identity component $\mbox{\it  Diff}_0(M)$
of the diffeomorphism group.  
A  beautiful  result of Donaldson  \cite[p. 336]{donperiod} asserts  that the
 period map  is  not only smooth, but is  actually transverse to the set of 
  planes containing any given element of positive self-intersection.
 This has  the following useful consequence: 
 
 \begin{prop}[Donaldson] \label{period} 
 Let $(M,g)$ be  any smooth compact oriented $4$-manifold with $b_+(M) \geq 1$,
 and let ${\mathbf b}\in {\mathcal H}^+_g\subset H^2(M,\RR )$ be the de Rham 
 class of any non-zero harmonic self-dual $2$-form on $(M,g)$. Then
 there is a smooth family of Riemannian metrics 
$g_{\mathbf t}$, ${\mathbf t} \in B_{\varepsilon} ({\mathbf 0})\subset {\mathcal H}^-_g$, 
with $g_{\mathbf 0}=g$,  such that 
$({\mathbf b}+{\mathbf t})\in {\mathcal H}^+_{g_{\mathbf t}}$
for each and every  ${\mathbf t}$. 
 \end{prop}

As the above discussion makes clear,  
the Hodge Laplacian
$$\Delta_d = dd^*+d^*d = - \star d \star d - d \star d \star$$
is an operator of fundamental geometric importance. It is thus worth pointing out that 
 if $\psi$ is a self-dual $2$-form, then  $\Delta_d\psi$ is also self-dual, and 
can, moreover, 
 be re-expressed by means of  the  Weitzenb\"ock formula \cite{bourg}   \begin{equation}
\label{friend}
\Delta_d \psi =     \nabla^{*}\nabla \psi - 2W_{+}(\psi , 
\cdot ) + \frac{s}{3} \psi  ~. 
\end{equation}
Taking the $L^2$ inner product with $\psi$, we therefore have 
$$
  \int_{M} \left(|\nabla \psi |^{2}
-2W_{+}(\psi , \psi ) + \frac{s}{3}|\psi |^{2}\right) d\mu\geq 0 ,  
$$
since $\Delta_d$ is a  non-negative operator. 
On the other hand, since $W_+: \Lambda^+\to \Lambda^+$ is self-adjoint and trace-free, 
$$|W_+(\psi , \psi ) | \leq \sqrt{\frac{2}{3}}|W_+| |\psi |^2,$$
so it follows that any self-dual $2$-form $\psi$ satisfies  
\begin{equation}
\label{part} 
 \int_{M} |\nabla \psi |^{2}d\mu \geq 
\int_M\left(-2\sqrt{\frac{2}{3}}|W_+|-\frac{s}{3}\right)|\psi |^{2} d\mu .
\end{equation}
Moreover,  
 assuming that $\psi\not\equiv  0$, 
equality holds iff $\psi$ is closed, belongs the lowest eigenspace of 
$W_+$ at each point, and the two largest eigenvalues of $W_+$ are
everywhere  equal. Of course, this last assertion crucially depends on  
the fact \cite{arons,baer} that if $\Delta_d \psi =0$ and $\psi \not\equiv 0$, then $\psi\neq 0$
on a dense subset of $M$.

A rather special set of techniques can be applied when $(M,g)$ happens to 
admit a closed self-dual 2-form $\omega\in {\mathcal H}^+_g$ with constant 
point-wise norm  $|\omega|_g \equiv \sqrt{2}$. In this case, there is 
an almost-complex structure $J:TM\to TM$, $J^2 =1$, defined by 
$$ g (J\cdot, \cdot ) = \omega (\cdot , \cdot ) ,$$
and this almost-complex structure then acts on $TM$ in a $g$-preserving fashion. 
The triple $(M,g,\omega )$ is then said to be an {\em almost-K\"ahler $4$-manifold}. 
Because $J$ allows one to to think of $TM$ as a complex vector bundle, 
it is only natural to look for a connection on its anti-canonical line bundle 
$L=\wedge^2 T^{1,0}_J\cong \Lambda^{0,2}_J$ in order to use 
the Chern-Weil theorem in order to  express $c_1^\RR (M,J)$ as 
$$ c_1^\RR (L)= [\frac{i}{2\pi}F]\in H^2_{DR}(M, \RR)~,$$ 
where $F$ is the curvature of
the relevant connection on $L$. A particular choice of 
Hermitian connection on $L$ was  first introduced by Blair \cite{blair}, and  is 
so geometrically natural that it was later rediscovered
by Taubes \cite{taubes} for entirely different reasons. The curvature
$F_B=F_B^++F_B^-$ of 
this {\em Blair connection}  is given  \cite{td,lsymp}  by 
\begin{eqnarray}
iF_B^+ &=& \frac{s+s^*}{8}\omega+  W^+(\omega )^\perp   \label{bla1}\\
iF_B^-&=&  \frac{s-s^\star }{8}\hat{\omega}+ \mathring{\varrho} \label{bla2}
\end{eqnarray}
where  the so-called {\em star-scalar curvature} is given by 
$$s^{*}=s+ |\nabla \omega |^2 = 2W_+ (\omega , \omega ) +  \frac{s}{3}~,$$
while   $W^+(\omega )^\perp$
  is the component of $W^+(\omega )$ orthogonal to $\omega$,  
  $$\mathring{\varrho}(\cdot , J\cdot )  = \frac{\mathring{r} +  J^*\mathring{r} }{2} ,$$
and where  the  anti-self-dual $2$-form 
 $\hat{\omega}\in \Lambda^-$
 is  defined only on the open set where $s^\star - s\neq 0$,  and   satisfies 
 $|\hat{\omega}|\equiv \sqrt{2}.$

An important special case occurs when $\nabla \omega=0$. This happens precisely
when 
$J$ is integrable, and $g$ is a  {\em K\"ahler metric} on the complex
surface $(M,J)$. In this case, $s=s^*$, $\omega$ is an eigenvector of the 
$W_+$, and $r$ is $J$-invariant, so that $iF_B$ just becomes the 
Ricci form $\rho$ defined by $\rho (\cdot , \cdot ) = r(J\cdot , \cdot )$. In fact,
$\omega$ is an eigenvector of $W_+$ with eigenvalue $s/6$, whereas the 
elements of $\omega^\perp= \Re e \Lambda^{2,0}_J$ are eigenvectors of eigenvalue $-s/12$. 

K\"ahler manifolds with  scalar curvature $s = \text{const} < 0$ will play 
an important r\^ole in this paper. By the above discussion, they belong to  
the following broader class of almost-K\"ahler
manifolds:

\begin{defn} \label{fatso} 
An almost-K\"ahler $4$-manifold $(M^4, g, \omega )$ will be said to be 
{\em saturated} if 
 \begin{itemize}
  \item $s+s^*$ is a negative constant; 
 \item $\omega$ belongs to the lowest eigenspace of $W_+: \Lambda^+\to \Lambda^+$ 
everywhere; and 
 \item the two largest eigenvalues of  $W_+: \Lambda^+\to \Lambda^+$ 
 are   everywhere equal.
 \end{itemize}
\end{defn}

\section{The Seiberg-Witten Equations} 
\label{monopoly}

This section is intended  both to fix 
our  terminological  conventions and to provide streamlined 
proofs of  the key preliminary
 curvature estimates. We note that, while all the main results 
in 
 this section  can be found elsewhere  \cite{il2,lno,lric,lebsurv2}, several
  of the proofs
given here are considerably simpler than those published heretofore. 

We begin with a discussion  of  {\em spin$^c$ structures}. 
If $M$ is any smooth  oriented $4$-manifold,  
 its second Stieffel-Whitney class $w_2(TM)\in H^2(M,\ZZ_2)$ is always \cite{hiho,reespinc}
 in the 
image of the natural homomorphism $$H^2 (M , \ZZ) \to H^2 (M , \ZZ_2),$$
and we can therefore always find
Hermitian line bundles $L\to M$ such that 
$$c_{1}(L)\equiv w_{2}(TM)
\bmod 2.$$
For any such $L$, and for any Riemannian metric $g$ on $M$,
one can then find rank-$2$ Hermitian vector bundles ${\mathbb V}_\pm$
which formally satisfy
$${\mathbb V}_\pm= {\mathbb S}_{\pm}\otimes L^{1/2},$$
where ${\mathbb S}_{\pm}$ are the locally defined left- and
right-handed spinor bundles of $(M,g)$. 
Such a choice of   ${\mathbb V}_{\pm}$, up to isomorphism, 
is called a 
spin$^{c}$ structure $\mathfrak c$ on $M$.
Moreover, if 
 $H_{1}(M,\ZZ)$ does not  contain elements of order $2$, then 
 $\mathfrak c$ is completely determined by 
$$c_{1}(L)= c_{1}({\mathbb V}_{\pm})
\in H^{2}(M,\ZZ ), $$ 
which is called the first Chern class of the spin$^c$ structure ${\mathfrak c}$. 

Every 
unitary 
connection $A$ on $L$ induces a connection 
$$\nabla_A : \Gamma ({\mathbb V}_{+})\to \Gamma (\Lambda^1\otimes {\mathbb V}_{+}),$$
and composition of this with the natural {\em Clifford multiplication} homomorphism
$$\Lambda^1\otimes {\mathbb V}_{+}\to {\mathbb V}_{-}$$
gives one a 
spin$^c$ version 
$$D_{A}: \Gamma ({\mathbb V}_{+})\to \Gamma ({\mathbb V}_{-})$$
of the Dirac operator \cite{hitharm,lawmic}. This is an elliptic first-order differential operator, 
and in many respects  it closely resembles the usual Dirac operator of 
spin geometry. In particular,  one has  the so-called Weitzenb\"ock formula 
\begin{equation}
\label{wtw} 
\langle \Phi , D_A^*D_A \Phi  \rangle = \frac{1}{2}\Delta |\Phi |^2 + |\nabla_A \Phi |^2 + 
\frac{s}{4} |\Phi |^2 + 2 \langle -iF_A^+ , \sigma (\Phi ) \rangle 
\end{equation}
for any  $\Phi\in \Gamma ({\mathbb V}_+)$, 
 where $F_A^+$ is the self-dual part of the 
curvature  of $A$,  and where 
$\sigma : {\mathbb V}_+ \to \Lambda^+$
is a natural real-quadratic map satisfying 
$$|\sigma (\Phi ) | = \frac{1}{2\sqrt{2}}|\Phi |^{2}.$$

Equation \eqref{wtw} is a natural  generalization of  the Weitzenb\"ock formula 
used by  Lichnerowicz \cite{lic} to prove that metrics with $s > 0$ cannot exist
when $M$ is spin and $\tau (M) \neq 0$. Unfortunately, however, one cannot hope to 
derive interesting  geometric information about the Riemannian
metric $g$ by just using \eqref{wtw} for an arbitrary connection $A$,
since one would have no  control at all over the  
 $F_A^+$ term.
 Witten \cite{witten}, however,  had the brilliant insight that 
one could remedy this by  considering both $\Phi$ and $A$
as unknowns, subject to the  
{\em Seiberg-Witten equations}  
\begin{eqnarray} D_{A}\Phi &=&0\label{drc}\\
-i F_{A}^+&=& \sigma(\Phi) .\label{sd}\end{eqnarray}
These equations are non-linear, but they become an 
 elliptic first-order system 
once one imposes the `gauge-fixing' condition
\begin{equation}
\label{gawk} 
d^* (A-A_0)=0,
\end{equation}
where $A_0$ is an arbitrary  `background''
connection on $L$, and 
$i(A-A_0)$ is simply treated as a real-valued 
$1$-form on $M$. 
The 
eliminate the natural action of the `gauge group' 
of  automorphisms of the Hermitian line bundle 
$L\to M$.

Because the Seiberg-Witten equations are non-linear, one cannot 
use something like an index formula to  predict that they
must have solutions. Nonetheless,  there exist  spin$^{c}$ structures on 
many  $4$-manifolds  for which 
there is at least one solution 
for every metric $g$. This situation is conveniently described  by the 
following 
 terminology \cite{K}:
 
\begin{defn} \label{kronos} 
Let $M$ be a smooth compact oriented $4$-manifold
with $b_{+}\geq 2$. An element ${\mathbf a}\in  H^{2}(M,\RR )$ is called a {\bf monopole
class} of $M$ iff there exists a spin$^{c}$ structure
$\mathfrak c$ 
on $M$ with  
$$c_{1}^\RR(L)=  {\mathbf a} $$ for which   the  Seiberg-Witten 
equations (\ref{drc}--\ref{sd})
have a solution for every Riemannian  metric $g$ on $M$. 
\end{defn}

When the gauge-fixing condition \eqref{gawk} is imposed, the 
Seiberg-Witten equations amount to saying that 
$(\Phi , A)$ belongs to the pre-image of zero for a Fredholm map of 
Banach spaces. This so-called {\em monopole map} turns
out to behave roughly like a proper map of finite-dimensional 
spaces \cite{baufu}. When the `expected dimension' of the
moduli space of solutions modulo gauge equivalence, as 
determined by the Fredholm index of the monopole map,
is zero,  
Witten \cite{witten}  discovered that 
one can define an invariant analogous to the degree of 
a map between finite-dimensional manifolds of the same dimension. 
More recently, Bauer and Furuta  \cite{baufu,bauer2}
discovered that the monopole map determines a
well-defined class in an equivariant cohomotopy group.
Either of these invariants can be used \cite{il2} to detect the presence 
of a monople class. Moreover, these invariants are often non-trivial; 
for example, a celebrated result of Taubes \cite{taubes} shows that 
if $(M,\omega )$ is a symplectic $4$-manifold with $b_+\geq 2$,  then 
Witten's invariant is non-zero for the spin$^c$ structure canoncially
detemined by $\omega$, so that 
$\pm c_1 (M,\omega)$ are both monopole classes. 
 On the other hand, Kronheimer \cite{K} has has used the Floer homology 
 of $3$-manifolds to show that 
 some $4$-manifolds admit monopole classes which are
 not detected by these invariants.

Equations (\ref{drc}--\ref{sd}) are precisely chosen  so as to  imply the 
 {Weitzenb\"ock formula}
 \begin{equation}
 0=	2\Delta |\Phi|^2 + 4|\nabla_A \Phi|^2 +s|\Phi|^2 + |\Phi|^4 .	
 	\label{wnbk}
 \end{equation}
In particular,
 these Seiberg-Witten  equations can never admit a solution 
 $(\Phi , A)$ with  $\Phi \not\equiv 0$
relative to a metric $g$ with  $s >0$. 
This leads, in particular, 
to a cornucopia of  
simply connected non-spin $4$-manifolds which do not admit 
positive-scalar-curvature metrics --- in complete contrast to the situation in 
higher dimensions \cite{gvln}. 
Even more strikingly,  the Seiberg-Witten  equations actually 
imply integral estimates for  the scalar curvature
\cite{witten,lpm}:

\begin{prop}\label{best}
Let $(M,g)$ be a smooth compact oriented Riemannian manifold,
and let ${\mathbf a}\in H^2(M, \RR)$ be a monopole class of $M$. 
Then the scalar curvature 
$s$ of $g$ satisfies 
$$
\int_{M}s^{2}d\mu_{g} \geq 32\pi^{2} ({\mathbf a}^{+})^{2} .
$$
If  ${\mathbf a}^+\neq 0$, moreover,  equality  occurs iff 
there is an integrable complex structure $J$ with 
 $c_1^\RR(M,J)={\mathbf a}$ such that $(M,g,J)$ is a 
 K\"ahler manifold of  constant negative scalar curvature. 
\end{prop}
\begin{proof}
By Definition \ref{kronos}, there must be a spin$^c$ structure with $c_1^\RR(L)={\mathbf a}$ 
for which the Seiberg-Witen equations (\ref{drc}--\ref{sd}) have
 a solution $(\Phi , A)$ on $(M,g)$. However,  given such a solution, 
$\Phi$  satisfies the Weitzenb\"ock formula 
 \eqref{wnbk} with respect to $g$ and $A$, and    integrating this  then reveals that 
$$0= \int [ 4|\nabla_A \Phi |^2 + s|\Phi|^2 + |\Phi|^4 ] d\mu . $$
Hence 
$$\int (-s) |\Phi|^2 d\mu \geq \int |\Phi|^4 d\mu ,$$
and applying the Cauchy-Schwarz inequality to the
left-hand side  yields 
$$
\left(\int s^2 d\mu  \right)^{1/2}\left(\int |\Phi |^4 d\mu  \right)^{1/2} \geq \int |\Phi |^4 d\mu . 
$$
Equation \eqref{sd} therefore tells us that 
$$
\int s^2 d\mu \geq \int |\Phi |^4 d\mu  =  8 \int |F_A^+|^2 d\mu . 
$$ 
However, since $iF_A/ 2\pi$ represents  ${\mathbf a}$ in de Rham cohomology, 
 Lemma \ref{harmless}
tells us that 
$$\int |F_A^+|^2 d\mu \geq 4\pi^2  ({\mathbf a}^+)^2~. $$
It follows  that 
$$
\int s^2 d\mu \geq 32\pi^2  ({\mathbf a}^+)^2~, 
$$
exactly as claimed.

If equality holds, the inequalities in the above argument must all be equalities.
Hence  $\nabla_A \Phi=0$, 
and   $iF_A^+ = -\sigma (\Phi )$  is therefore  a parallel self-dual $2$-form with 
de Rham class  $2\pi {\mathbf a}^+$. When this cohomology class is non-zero,
this form cannot vanish, and 
we therefore conclude that $(M,g)$ is K\"ahler. Inspection of \eqref{wnbk}
then reveals  that  $s$ must be a negative constant. 
Moreover, $\Phi \otimes \Phi$ is then a non-zero section of $\Lambda^{2,0}\otimes L$
with respect to the relevant  complex structure, so $c_1^\RR (M, J) = c_1^\RR (L) = {\mathbf a}$. 
Conversely,  any 
such K\"ahler metric  would saturate the inequality because  
the self-dual part of the Ricci form of any K\"ahler metric on a K\"ahler surface is 
 $s \omega /4$, where 
$\omega = g (J\cdot , \cdot )$ is the K\"ahler form. 
\end{proof}

\begin{prop}
\label{sneg} 
Let $M$ be a compact oriented $4$-manifold with $b_+(M) \geq 2$. 
If there is a non-zero monopole class ${\mathbf a}\in H^2 (M,  \RR)-\{ {\mathbf 0}\}$, 
 then $M$ does not admit  metrics
of scalar curvature $s\geq 0$. 
\end{prop}

\begin{proof}
Let $M$ be a smooth compact $4$-manifold with $b_+(M) \geq 2$, and 
suppose that ${\mathbf a} \in H^2(M, \RR)-{\mathbf 0}$
is a non-zero monopole class. By definition, this means that here is a 
spin$^c$ structure on $M$ with $c_1^\RR(L)= {\mathbf a}$
for which the Seiberg-Witten equations  have some  solution 
$(\Phi , A)$ with respect to any metric $g$ on $M$. But if the metric in question has
$s\geq 0$, the Weitzenb\"ock formula \eqref{wnbk}  says that
$$0 = 2\Delta |\Phi|^2 + |\nabla_A \Phi |^2 + s |\Phi |^2 + |\Phi |^4$$
so that $s\geq 0$ implies that 
$$0\geq \int |\Phi|^4 d\mu_g$$
and we therefore have $\Phi \equiv 0$. But this implies that 
$F_A^+ = i \sigma (\Phi )\equiv 0$, too, so that ${\mathbf a} = [ \frac{i}{2\pi} F_A] \in {\mathcal H}^-_g$.
In particular, if $g$ has scalar curvature $s\geq 0$ and if ${\mathbf b}\in {\mathcal H}^+_g$,
then ${\mathbf a}\cdot {\mathbf b} =0$.

Next, suppose that we had some  metric $g$ on $M$  with strictly positive scalar curvature 
$s> 0$. 
Choose some  ${\mathbf b}\in {\mathcal H}^+_g$ with ${\mathbf b}^2 =1$. 
The argument in the previous paragraph tells us that ${\mathbf a} \in {\mathcal H}^-_g$, 
so that the integer class ${\mathbf a}\neq {\mathbf 0}$  must satisfy ${\mathbf a}^2 \leq -1$. 
However, Proposition \ref{period} now tells us we  can now find a smooth $1$-parameter family
of metrics $g_t$, $t\in (-\varepsilon , \varepsilon )$, 
 such that $g_0=g$, and such that ${\mathbf b} + t {\mathbf a}\in {\mathcal H}^+_{g_t}$
 for all $t$. 
Since we have assumed that $g$ has $s > 0$, the same is necessarily true
of all the metrics  $g_t$ for  sufficiently small $t$, and we thus 
certainly have a contradiction, 
 since the argument of the previous paragraph would now tell us that 
 ${\mathbf a} \cdot ( {\mathbf b }+t {\mathbf a}) = t {\mathbf a}^2$
 would have to vanish for all small values of $t$. 
 It follows that $M$ cannot admit any metrics of positive scalar curvature. 
 
Finally, let us suppose  instead that   $g$ is  a metric on $M$ with $s\geq 0$. 
Since  $M$ is now known  not  to admit  metrics of positive scalar curvature, 
$g$ must then have $s\equiv 0$, since otherwise \cite{aubis,bes} we would be able to 
produce a metric of strictly positive scalar curvature by conformally rescaling it. 
 We may now proceed much as in the previous case. 
Once again,  $s\equiv 0$ implies that 
 ${\mathbf a} \in {\mathcal H}^-_g$. Again, choose some  
 ${\mathbf b}\in {\mathcal H}^+_g$ with ${\mathbf b}^2 =1$,
and observe that, once again,
 there exists a family of metrics $g_t$, $t\in (-\varepsilon , \varepsilon )$ 
with $g_0=g$ for which  ${\mathbf b} + t {\mathbf a}\in {\mathcal H}^+_{g_t}$.
But this time, we invoke  a 
theorem of Koiso \cite{bes,koiso}  on the Yamabe problem with parameters, and thereby  
   construct a smooth family 
 of constant-scalar-curvature, unit-volume metrics $\tilde{g}_t$ by 
 conformally rescaling each $g_t$. The conformal invariance of 
 \eqref{deco} then tells us that we still have 
    ${\mathbf b} + t {\mathbf a}\in {\mathcal H}^+_{\tilde{g}_t}$. Since 
   the family $\tilde{g}_t$ is smooth, the value $s_{\tilde{g}_t}$
   of its scalar curvature is therefore a smooth function 
   of $t$. But since $M$ does not admit metrics of positive scalar curvature,
   and since $s_{\tilde{g}_0}=0$, 
   this smooth function must have a maximum at $t=0$. Hence there is a positive constant $C$ 
   such that 
   $$
   0\geq s_{\tilde{g}_t} \geq - Ct^2
   $$
  for all sufficiently small $t$, and we therefore have 
   $$C t^4\geq   s_{\tilde{g}_t}^2= \int_M  s_{\tilde{g}_t}^2 d\mu_{\tilde{g}_t}  $$
   for $t$ in the same range. However, Proposition \ref{best} and the 
   Cauchy-Schwarz inequality tell us that 
   $$
  \int_M  s_{\tilde{g}_t}^2 d\mu_{\tilde{g}_t} \geq 32\pi^2 ({\mathbf a}^+_{\tilde{g}_t})^2
  \geq  32\pi^2    \frac{[{\mathbf a}\cdot ({\mathbf b}+ t{\mathbf a})]^2}{({\mathbf b}+ t{\mathbf a})^2}
  = 32\pi^2  \frac{ t^2 |{\mathbf a}^2|^2}{1 -t^2 |{\mathbf a}^2|}\geq 32\pi^2  t^2
   $$
   so we conclude that  $(\text{const})  t^4\geq t^2$ for all small $t$, which is certainly a 
   contradiction. Thus no metric with $s\geq 0$ can exist, and we are done.  
 \end{proof}

\begin{defn}
For any   smooth compact oriented $4$-manifold $M$ 
with $b_{+}\geq 2$, we  set  
$${\mathfrak C}(M) = \left\{ \text{monopole classes }  {\mathbf a} \in H^2(M, \RR)
\right\}. $$
We will often  abbreviate ${\mathfrak C}(M)$ as ${\mathfrak C}$ when no confusion seems likely to result. 
\end{defn}

\begin{lem}\label{libra} 
For   any smooth compact oriented $4$-manifold $M$ with $b_{+}\geq 2$,
$${\mathfrak C}(M)  = - {\mathfrak C}(M).$$
That is,   ${\mathbf a}\in  H^{2}(M,\RR )$ is a  monopole
class  iff  $-{\mathbf a}\in  H^{2}(M,\RR )$ is a  monopole class, too. 
\end{lem} 
\begin{proof}
Let $g$ be any metric  on $M$, and let  ${\mathbb V}_\pm$
be the twisted spin bundles of some spin$^c$ structure $\mathfrak c$. Then
the conjugate vector bundles $\overline{\mathbb V}_\pm$ are the 
twisted spin bundles of a second spin$^c$ structure $\overline{\mathfrak c}$,
since we have  natural isomorphisms 
$$
\overline{\mathbb V}_\pm\cong {\mathbb V}_\pm \otimes L^{-1} , ~~~~~L=\det ({\mathbb V}_\pm)
$$
induced by the wedge and inner products. Since we locally  have
\begin{eqnarray*}
{\mathbb V}_\pm&=&  {\mathbb S}_\pm\otimes L^{1/2}
\\
\overline{\mathbb V}_\pm&=&  {\mathbb S}_\pm\otimes L^{-1/2}
\end{eqnarray*}
as bundles with connection, it follows that 
$$
\overline{D_A \Phi} = D_{\bar{A}} \overline{\Phi}
$$
for any $\Phi \in \Gamma ({\mathbb V}_+)$ and any 
Hermitian connection $A$ on $L$, where $\bar{A}$ denotes the dual connection on 
$L^{-1}$ induced by $A$. 
Moreover, since the associated anti-linear map 
$$
{\mathbb S}_+\to {\mathbb S}_+
$$
acts by multiplying by the quaternion $j$, we also have 
$$
\sigma (\overline{\Phi})= - \sigma (\Phi ) .
$$
Since $F_{A*}=-F_A$, it follows that 
 if $(\Phi , A)$ is a solution of (\ref{drc}--\ref{sd})  with respect to 
$(g, {\mathfrak c})$,  then $(\overline{\Phi} , \bar{A})$ 
is a solution of (\ref{drc}--\ref{sd}) with respect to  $(g, \overline{\mathfrak c})$. 
If, for every metric $g$ on $M$,  there is a solution of the Seiberg-Witten equations 
for the spin$^c$ structure ${\mathfrak c}$, the same is therefore also true of  $\overline{\mathfrak c}$.
Since 
$c_1(\overline{\mathfrak c}) = c_1(\overline{\mathbb V}_+) = - c_1({\mathbb V}_+)
=-c_1({\mathfrak c})$, it follows that 
the set of monopole classes is invariant under multiplication 
by $-1$. 
\end{proof}

A particularly important consequence  of Proposition \ref{best}  is the following fundamental 
fact \cite{il2}:

\begin{prop}
 Let $M$ be any smooth compact oriented $4$-manifold with 
 $b_+(M)\geq 2$. Then  ${\mathfrak C}(M )$ 
  is a finite set. 
 \end{prop}
 \begin{proof} First, observe that one can always  find a basis 
 $\{ e_j~| j=1, \ldots , b_2\}$ for  $H^2(M,\RR )$, together with a  collection of 
 metrics $g_j$ such that $e_j\in {\mathcal H}^+_{g_j}$. 
Indeed, let $g= g_1= \ldots = g_{b_+}$ be any metric, let
 $e_1, \ldots , e_{b_+}$ to be any basis for
 ${\mathcal H}^+_{g}$, and let  $e_{b_++1}, \ldots , e_{b_2}$ then 
 be small perturbations of $e_1$ by linearly independent elements of ${\mathcal H}^-_{g}$,
 while using Proposition \ref{period} to find compatible metrics 
 $g_{b_++1}, \ldots , g_{b_2}$. 
 Alternatively, 
 one can simply  take the $e_j$ to be any collection of rational classes 
with $e_j^2 >0$  which span  $H^2(M,\RR )$, and then cite 
 a remarkable recent  construction of Gay and Kirby  \cite[Theorem 1]{kirbyperiod}, 
 which 
 shows that any rational cohomology class with positive self-intersection can be
 be represented by a closed $2$-form which is self-dual with respect to some metric.
 Given this data, 
 we now introduce a   constant  for each $j$ by setting 
$$
\kappa_j = \left(\frac{e_j^2}{32\pi^2}\int_M s_{g_j}^2d\mu_{g_j}\right)^{1/2}. 
$$

 Let $L_j: H^2(M,\RR ) \to \RR$
be the linear functionals $L_j(x)={e}_j\cdot  x$. 
Since the intersection form is positive-definite on  
each 
${\mathcal H}^+_{g_j}$, 
the Cauchy-Schwarz inequality and 
Proposition 
\ref{best} together
imply   that any monopole class 
${\mathbf a}\in H^2(M, \RR )$ must satisfy
$$\left|L_j ({\mathbf a})\right| = \left| e_j \cdot {\mathbf a}\right|
=
\left| e_j\cdot {\mathbf a}^+_{g_j}\right|\leq \sqrt{e_j^2}\sqrt{({\mathbf a}^+)^2}
\leq \kappa_j$$
for each $j$. 
Hence 
 ${\mathfrak  C}\subset H^2(M,\RR )$
is contained in the parallelepiped 
$$\left\{ x\in H^2 (M,\RR ) ~\Big|~ |L_j(x)|\leq \kappa_j ~\forall j= 1,\ldots , b_2(M)\right\},$$
which is  a compact set. 
But ${\mathfrak C}\subset H^2(M,\ZZ )/\mbox{torsion}$, and  is
therefore also discrete. Hence $\mathfrak C$ is finite, as claimed. 
 \end{proof}
 
We now introduce a generalization of the  Seiberg-Witten equations.
Let $(M,g)$  be a smooth oriented Riemannian $4$-manifold, let
${\mathfrak c}$ be a spin$^c$-structure on $M$, and let $f: M\to \RR^+$ be a
smooth positive function. Then
 we will say that  $(\Phi , A)$ solves the {\em rescaled Seiberg-Witten equations}
if 
\begin{eqnarray} D_{A}\Phi &=&0\label{rsdrc}\\
-i F_{A}^+&=&f\sigma (\Phi) \label{rssd}\end{eqnarray}

\begin{lem} \label{hilfsatz}
Let $M$ be a smooth compact $4$-manifold with $b_+\geq 2$, and let 
${\mathbf a}\in H^2(M, \RR)$ be a monopole class. Then, for any smooth metric $g$ and
any smooth positive function $f$, there is a solution of the rescaled Seiberg-Witten 
equations (\ref{rsdrc}--\ref{rssd}) for some spin$^c$ structure on 
$M$ with $c_1^\RR (L) = {\mathbf a}$. 
\end{lem}
\begin{proof}
Consider the conformally related metric $\hat{g}=f^{-2}g$. 
Because ${\mathbf a}$ is a monopole class,  there must then be a 
solution $(\hat{\Phi}, A)$ of the Seiberg-Witten equations wtih respect to  $\hat{g}$ and some 
spin$^c$ structure with $c_1^\RR (L) = {\mathbf a}$. However, the Dirac equation \eqref{drc}
is conformally invariant. More precisely, $\hat{\Phi}$ uniquely determines \cite{lawmic,pr1}
 a solution 
$\Phi$ of \eqref{drc} with respect to $g$, such that $|\Phi|_g = f^{-3/2}|\hat{\Phi}|_{\hat{g}}$, 
and such that $\sigma_g (\Phi) = f^{-1} \sigma_{\hat{g}}(\hat{\Phi})$.  Hence
$(\Phi , A)$ satisfies (\ref{rsdrc}--\ref{rssd}) with respect to $g$. 
\end{proof}

Given a solution  $(\hat{\Phi}, A)$ of (\ref{rsdrc}--\ref{rssd}), substitution 
into  (\ref{wtw}) yields   the Weitzenb\"ock formula 
 $$
0 = 2\Delta |\Phi |^{2} + 4 |\nabla_{A} \Phi |^{2} + s|\Phi |^{2}
 + f|\Phi |^{4} .
$$
Multiplying by $|\Phi |^2$ and integrating, we then obtain 
an  inequality 
\begin{equation}
\label{two}
0\geq \int_M\left[4 |\Phi |^{2}|\nabla_{A} \Phi |^{2} + s|\Phi |^{4}
 + f|\Phi |^{6}\right] d\mu 
\end{equation}
and we will now use this to prove our next result.

 \begin{prop}\label{upscale}
 Let $(M,g)$ be a smooth compact oriented Riemannian manifold,
 and let ${\mathbf a}\in H^2(M, \RR)$ be a monopole class of $M$. 
Then the scalar curvature 
$s$ and self-dual Weyl curvature $W_+$ of $g$ satisfy
 $$\int_M (s-\sqrt{6}|W_+|)^2d\mu_g \geq 72\pi^2({\mathbf a}^+)^2.$$
 If  ${\mathbf a}^+\neq 0$,  moreover, equality holds iff 
 there is a   symplectic form $\omega$, where $[\omega ]$ is a negative multiple of ${\mathbf a}^+$ 
 and   $ c_1^\RR (M , \omega )={\mathbf a}$, such that
 $(M,g,\omega )$ is a saturated 
 almost-K\"ahler  manifold in the sense of Definition \ref{fatso}, 
 . 
  \end{prop}
 \begin{proof} 
For any   smooth  function $f> 0$ on $M$, Lemma \ref{hilfsatz} 
 guarantees that the 
corresponding rescaled 
Seiberg-Witten equations (\ref{rsdrc}--\ref{rssd}) must 
have some solution  $(\Phi ,  A)$.  Set   $\psi = 2\sqrt{2} \sigma (\Phi )$, and 
 observe that   the definition of $\sigma$ then implies that 
$$
|\Phi |^{4} =	|\psi |^{2}   ,~~~~~~~~
	   4 |\Phi |^{2}|\nabla_A \Phi |^{2}  \geq  |\nabla \psi |^{2}  . 
$$
Thus inequality  (\ref{two})
 tells us that 
 $$0\geq \int_M\left[|\nabla\psi  |^{2} + s|\psi |^{2}
 + f|\psi |^{3}\right] d\mu . $$
However, inequality  (\ref{part}) also tells us that
$$
\int_{M} |\nabla \psi |^{2}d\mu \geq 
\int_M\left(-2\sqrt{\frac{2}{3}}|W_+|-\frac{s}{3}\right)|\psi |^{2} d\mu , 
$$
and combining these facts yields 
$$
0\geq \int_M\left[  \left(\frac{2}{3} s-2\sqrt{\frac{2}{3}}|W_+|\right) |\psi |^{2} 
 + f|\psi |^{3}\right] d\mu .
$$
Set $\varphi = \frac{3}{2}\psi = 3\sqrt{2}\sigma (\Phi )$. We then have
$$
0\geq \int_M\left[\left( s-\sqrt{6}|W_+|\right) |\varphi |^{2} 
 +   f|\varphi |^{3}\right] d\mu .
$$
Rewriting this as 
$$
\int_M\left[ -\left( s-\sqrt{6}|W_+|\right)f^{-2/3}\right] \left( f^{2/3}|\varphi  |^{2}\right) d\mu \geq 
\int_M f|\varphi  |^{3}  d\mu 
$$
and applying the H\"older inequality to the left-hand side then yields 
    $$
  \left[\int_M \left|s-\sqrt{6}|W_+|\right|^3f^{-2}d\mu\right]^{1/3}
  \left[ \int_M f|\varphi  |^{3}d\mu \right]^{2/3}
   d\mu \geq \int_M f|\varphi  |^{3}d\mu ,
   $$
  which is to say that 
    $$
  \int_M \left|s-\sqrt{6}|W_+|\right|^3f^{-2}d\mu   \geq 
   \int_M f|\varphi  |^{3}d\mu .
   $$
    But the H\"older inequality also tells us that 
    $$
   \left(\int_M f^4d\mu\right)^{1/3}   \left( 
   \int_M   f|\varphi  |^{3}d\mu\right)^{2/3} \geq 
  \int_M f^{4/3} \left[  f^{2/3} |\varphi  |^2\right] d\mu ~,    $$
   where equality holds only if $|\varphi |$ is a constant multiple of $f$. 
 Hence 
  $$
   \left(\int_M f^4d\mu\right)^{1/3}   \left( \int_M \left|s-\sqrt{6}|W_+|\right|^3
   f^{-2}d\mu\right)^{2/3} \geq 
  \int_M f^2 |\varphi |^2 d\mu ~.    $$
   But since $f\varphi = 3\sqrt{2} f\sigma (\Phi ) = 3\sqrt{2}(-iF_A^+)$,  we also  have 
   $$ \int_M f^2 |\varphi|^2 d\mu = 18 \int_M |F_A^+|^2 d\mu \geq 
  18 (2\pi {\mathbf a}^+)^2 
  =
   72\pi^2 ({\mathbf a}^+)^2 $$
   by Lemma \ref{harmless}, since $iF_A\in 2\pi c_1^\RR (L) = 2\pi {\mathbf a}$. 
   Thus 
   \begin{equation}\label{clem} 
\left(\int_M f^4d\mu_g \right)^{1/3}\left( \int_{M}\left|s-\sqrt{6}|W_+|\right|^{3}f^{-2}d\mu_g 
 \right)^{2/3}\geq  72\pi^2 ({\mathbf a}^+)^2 
\end{equation}
 for any smooth positive function $f$ on $M$.
 
 Now choose a sequence of smooth positive functions $f_j$ on 
 $M$ with 
 $$ f_j \searrow \sqrt{\left|  s-\sqrt{6}|W_+|   \right|}$$
uniformly on $M$.  Since the inequality $f_j^2 \geq {\left|  s-\sqrt{6}|W_+|   \right|}$
 implies 
$$
\int_M f_j^4 d\mu \geq  \left(\int_M f_j^4d\mu_g \right)^{1/3}\left( \int_{M}\left|s-\sqrt{6}|W_+|\right|^{3}f_j^{-2}d\mu_g 
 \right)^{2/3},  $$
we  then have 
$$
\int_M f_j^4 d\mu \geq  72\pi^2 ({\mathbf a}^+)^2 
$$
 by applying  (\ref{clem}). But since 
 $$
 \int_M \left(s-\sqrt{6}|W_+|\right)^{2}d\mu = \lim_{j\to \infty} 
 \int_M f_j^4 d\mu , 
 $$
this shows that 
\begin{equation}
\label{voila}
\int_M \left(s-\sqrt{6}|W_+|\right)^{2}d\mu  \geq  72\pi^2 ({\mathbf a}^+)^2 , 
\end{equation}
as  desired. 
 
Finally, we analyze the equality case. Suppose that $g$ is a metric such that 
equality holds in (\ref{voila}). Then $g$ must in particular minimize
$${\mathcal A}(g) = \int (s_g-\sqrt{6}|W_+|_g)^2d\mu_g$$ in its conformal class. However, 
if $u$ is any smooth positive function, and if $\hat{g}= u^2 g$, then
$$
{\mathcal A}(u^2g) = 
  \int (s_{{g}}+ 6 u^{-1}\Delta_g u -\sqrt{6}|W_+|_{{g}})^2 d\mu_{{g}}
$$
so that, for the $1$-parameter family of metrics given by 
$$g_t = (1+ t F)^2 g$$
one has 
$$
\frac{d}{dt}{\mathcal A}(g_t) \Big|_{t=0} = 
12 \int [\Delta_g F  ] (s_{\hat{g}}-\sqrt{6}|W_+|_{\hat{g}}) d\mu_{{g}} .
$$
If $g$ minimizes ${\mathcal A}$ in its conformal class, 
we must therefore have 
$$\Delta_g \left( s-\sqrt{6}|W_+| \right)=0$$
in the weak (or distributional) sense. Elliptic regularity  \cite{giltrud} therefore tells us that 
$s-\sqrt{6}|W_+|$ is smooth, and integrating by parts 
$$\int \left|\nabla  \left( s-\sqrt{6}|W_+| \right)\right|^2d\mu = 
\int  \left( s-\sqrt{6}|W_+| \right) \left[\Delta \left( s-\sqrt{6}|W_+| \right)\right]d\mu =0$$
therefore shows  that 
$$ s-\sqrt{6}|W_+| = \mbox{constant.}$$
Assuming ${\mathbf a}^+\neq 0$,  moreover, Proposition \ref{sneg} 
tells us this constant must be negative. 
With this proviso, we can then set 
$$f =  \sqrt{\left|  s-\sqrt{6}|W_+|   \right|} , $$
and equality in (\ref{voila}) then implies that equality occurs in (\ref{clem}) for this
choice of $f > 0$. But then, for this choice of $f$, we must therefore have 
equality at every step of the proof of 
 (\ref{clem}). Since this $f$ is constant, it thus follows that $\varphi= 3\sqrt{2}\sigma (\Phi )$ 
 is a closed self-dual $2$-form of non-zero constant length. 
 Setting
 $\omega =  \sqrt{2} \psi / |\psi |$, it   follows that 
  $(M, g,\omega)$
 is an   almost-K\"ahler manifold
 Moreover, since 
 $\psi=\frac{2}{3}\varphi$ belongs to the lowest eigenspace
 of $W_+$ at each point, while the two largest eigenvalues of $W_+$ must be equal
 at every point,
 we  have 
 $$|W_+| = \sqrt{\frac{3}{2}}\frac{[-W_+(\omega , \omega )]}{|\omega|^2} = -\frac{1}{2}  \sqrt{\frac{3}{2}}~
 W_+(\omega , \omega )$$
so that 
$$s + s^* = s+ \left[ \frac{s}{3} + 2W_+(\omega , \omega ) \right] = 
 \frac{4}{3} \left( s- \sqrt{6}|W_+|\right) , $$
which we already know to be  a negative constant. The almost-K\"ahler manifold  
$(M, g,\omega)$
 is  therefore saturated   in the sense of Definition \ref{fatso}. 
Moreover, since $\Phi \otimes \Phi$ is  a non-zero section of $\Lambda^{2,0}_J\otimes L$
we have $c_1^\RR (M, \omega) = c_1^\RR (L) = {\mathbf a}$. 
Moreover,  by construction, $\omega$ is a negative multiple of $iF_A^+/2\pi$, which is 
the harmonic representative of ${\mathbf a}^+$. 

Conversely, if $(M,g,\omega)$ is an almost-K\"ahler manifold with 
$b_+\geq 2$, then ${\mathbf a}= c_1^\RR (M, \omega )$ is a monopole
class by Taubes'  theorem \cite{taubes}, and in the saturated case 
our formula 
\eqref{bla1} then shows not only  that the harmonic representative of 
${\mathbf a}^+$ is given by $iF_B^+/{2\pi}$, where $F_B$ is the 
curvature of the Blair connection, but also moreover that equality 
occurs in \eqref{voila} for this choice of monopole class.
 The Proposition therefore follows. 
 \end{proof}
   
\section{Monopoles and Convex Hulls}

In the previous section, we saw that  monopole classes lead to 
non-trivial lower bounds for the $L^2$-norms of certain curvature 
expressions. Unfortunately, however, these lower bounds
still  depend on the image of $g$ under the period map, and so
are not yet  uniform in the metric. We will now remedy this,  using
 some simple tricks from  convex geometry.  

We begin by establishing a notational convention:

\begin{defn} 
Let ${\mathbb V}$ be a vector space over $\RR$, and 
let $S\subset {\mathbb V}$. Then $\hull (S)\subset {\mathbb V}$
will denote the {\em convex hull} of $S$, meaning the  
smallest convex subset of ${\mathbb V}$ which contains 
$S$. 
\end{defn}

\begin{lem} \label{reflect} 
Let  $M$ be a smooth compact oriented $4$-manifold with 
$b_+\geq 2$,  and let ${\mathfrak C}= {\mathfrak C}(M) \subset H^2(M, \RR)$
be its set of non-zero monopole classes. Then 
$\hull ({\mathfrak C})\subset H^2(M, \RR)$ is  compact. Moreover, 
$\hull ({\mathfrak C})$  is symmetric, in the sense that 
$\hull ({\mathfrak C})= -\hull ({\mathfrak C})$. 
\end{lem}
\begin{proof}
By definition, $\hull ({\mathfrak C})$ is the smallest convex subset of $H^2 (M, \RR)$
which contains ${\mathfrak C}(M)$. However, since ${\mathfrak C}(M)$
is a finite subset, say $\{  {\mathbf a}_1 , \ldots , {\mathbf  a}_n\}$, we can explicitly 
express this convex hull as 
$$
\hull ({\mathfrak C}) = \left\{ ~\sum_{j=1}^n t_j {\mathbf a}_j ~\Big|~
 t_j \in [0,1] , ~\sum_{j=1}^n t_j =1~\right\} ,  
$$
since the set on the right is certainly a convex subset containing the ${\mathbf a}_j$,
and conversely is necessarily contained in any convex subset containing these points. 
In particular, this means that $\hull ({\mathfrak C})$ is the image of the standard 
$(n-1)$-simplex 
$$\triangle^{n-1} = \left\{ (t_1 , \cdots ,t_n ) \in [0,1]^{n}
~ \Big| ~ \sum_j t_j = 1
 \right\} $$
 under the continuous map 
 $$
 (t_1 , \cdots ,t_n )\longmapsto \sum_{j=1}^n t_j {\mathbf a}_j ~,
 $$
 and, since $\triangle^{n-1}$ is compact, it follows that $\hull ({\mathfrak C})$
 is, too.
 
 On the other hand, Lemma  \ref{libra} tells us that ${\mathfrak C}(M)$ is symmetric.
 Hence 
 $$\hull ({\mathfrak C})=\hull (-{\mathfrak C})=-\hull ({\mathfrak C})$$
 and $\hull ({\mathfrak C})$ is therefore symmetric, too.
 \end{proof}

Let us now consider the self-intersection function
\begin{eqnarray*}
Q: H^2 (M, \RR ) &\longrightarrow& \RR\\
{\mathbf v} &\longmapsto&  {\mathbf v}^2 ~,
\end{eqnarray*}
where ${\mathbf v}^2$ is of course just short-hand for 
${\mathbf v}\cdot {\mathbf v}= \langle {\mathbf v}\smallsmile {\mathbf v}, [M]\rangle$. 
Notice that  $Q$ 
 is a polynomial function, and therefore  continuous. Since $\hull ({\mathfrak C})$ is compact
by Lemma \ref{reflect}, it thus follows that $Q|_{\hull ({\mathfrak C})}$ necessarily achieves its
 maximum.  
We are thus entitled to make the following definition: 

\begin{defn} \label{defbeta}
Let  $M$ be a smooth compact oriented $4$-manifold with 
$b_+\geq 2$, 
and let $\hull ({\mathfrak C})
\subset H^2 (M, \RR)$ denote the  convex hull
of the set ${\mathfrak C}= {\mathfrak C}(M)$ of monopole classes  of $M$.
If ${\mathfrak C} \neq {\varnothing}$, 
we define 
$$
{\beta}^2 (M) = \max \left\{ \left. {\mathbf v}^2 ~\right|~ {\mathbf v} \in \hull ({\mathfrak C}) ~\right\}.
$$
 If, on the other hand,  ${\mathfrak C}={\varnothing}$, we instead set ${\beta}^2 (M)=0$.
\end{defn}

\begin{prop} \label{yobo} 
For any smooth $M^4$ with $b_+\geq 2$, ${\beta}^2 (M) \geq 0$. 
\end{prop}
\begin{proof}
If ${\mathfrak C}={\varnothing}$,  we have ${\beta}^2 (M)=0$ by Definition \ref{defbeta}. 
Otherwise, let ${\mathbf a}\in {\mathfrak C}$, and observe that 
$-{\mathbf a}\in {\mathfrak C}$, too, by Lemma \ref{libra}. Thus 
${\mathbf 0}= \frac{1}{2}{\mathbf a} +\frac{1}{2} (-{\mathbf a}) \in \hull ({\mathfrak C})$. 
Hence 
$${\beta}^2 (M) = \max \left\{ \left. {\mathbf v}^2 ~\right|~ {\mathbf v} \in \hull ({\mathfrak C}) ~\right\}\geq
{\mathbf 0}^2 = 0,$$
exactly as claimed. 
\end{proof}

\begin{prop}\label{biga}
Let $M$ be a smooth compact oriented $4$-manifold with 
 ${\mathfrak C}(M)\neq \varnothing$. 
 Then, for any   Riemannian metric $g$ on $M$, there is
a monopole class ${\mathbf a}\in {\mathfrak C}(M)$ such that 
$$({\mathbf a}^+)^2 \geq \beta^2 (M).$$
\end{prop}
\begin{proof}
Let ${\mathbf v}\in \hull ({\mathfrak C})$ be a maximum point of 
$Q$, so  that 
 ${\mathbf v}^2=\beta^2 (M)$
by Definition \ref{defbeta}. 
Let $\Pi : H^2 (M, \RR) \to {\mathcal H}^+_g$ denote the orthogonal
projection map. Since $\Pi$ is a linear map, we automatically  have
 $\hull (\Pi ({\mathfrak C }))= \Pi (\hull ({\mathfrak C}))$. However, 
 since the intersection form is positive definite on 
  ${\mathcal H}^+_g$, $Q|_{{\mathcal H}^+_g}$  has positive-definite Hessian, and
  the maximum of $Q$ on a line segment $\overline{\mathbf p\mathbf q} \subset {\mathcal H}^+_g$ can 
  therefore never 
  occur at an  interior point. The maximum points  of $Q|_{\Pi (\hull ({\mathfrak C}))}$
 must therefore all belong to $\Pi ({\mathfrak C })$. In particular, there must be a 
  monopole class ${\mathbf a} \in {\mathfrak C}$ such that 
  $$({\mathbf a}^+)^2= Q(\Pi ({\mathbf a})) \geq Q(\Pi ({\mathbf v})) =({\mathbf v}^+)^2.$$
  On the other hand, 
  $${\mathbf v}^2 = ({\mathbf v}^+)^2 - |({\mathbf v}^-)^2|~,$$
  so we therefore have
  $$({\mathbf a}^+)^2\geq ({\mathbf v}^+)^2 \geq {\mathbf v}^2=  \beta^2 (M),$$
and the monopole class ${\mathbf a}$ therefore fulfills our desideratum. \end{proof}

The first part of Theorem \ref{summa} now follows immediately:

\begin{thm} \label{magna} 
Let $M$ be a smooth compact oriented $4$-manifold with $b_+\geq 2$. 
Then any metric $g$ on $M$ satisfies  curvature estimates \eqref{mocha} and \eqref{java}{\rm :} 
\begin{eqnarray*}
\int_M s^2d\mu &\geq&  32\pi^2 \beta^2(M) \\
\int_M (s-\sqrt{6}|W_+|)^2d\mu &\geq& 72\pi^2 \beta^2(M) 
\end{eqnarray*}
\end{thm}
\begin{proof}
For any metric $g$ on $M$, Proposition \ref{biga} tells us that 
there is a monopole class ${\mathbf a}$ such that 
$({\mathbf a}^+)^2 \geq \beta^2(M)$.  Proposition \ref{best} therefore tells us  that 
$$\int_M s^2d\mu  \geq  32\pi^2  ({\mathbf a}^+)^2 \geq  32\pi^2 \beta^2(M) ~,$$
while Proposition \ref{upscale} 
tells us that 
$$
\int_M (s-\sqrt{6}|W_+|)^2d\mu  \geq 72\pi^2  ({\mathbf a}^+)^2 \geq 72\pi^2 \beta^2(M) ,~
$$
and the Theorem therefore follows. 
\end{proof}

To prove Theorem \ref{summa}, we therefore merely need to understand the
equality cases of the curvature estimates  \eqref{mocha} and \eqref{java}. 
To do this, we will first need the following simple observation: 
 
 \begin{lem} \label{betapos} 
 Suppose that $(M,g)$ is a Riemannian manifold with 
 $b_+\geq 2$, and that $M$ carries a non-zero monopole class. 
  If equality occurs in either \eqref{mocha} or \eqref{java}, 
 then $\beta^2 (M) > 0$. 
 \end{lem}
 \begin{proof} If equality were to hold in \eqref{mocha} or \eqref{java}, and if we also had 
 $\beta^2(M) =0$, 
 the metric in question would necessarily have $s\geq 0$. But Proposition  \ref{sneg}
 says that no such metric can exist in the presence of a non-zero monopole class. 
 The claim thus follows by contradiction. 
 \end{proof}
 
 We will also need the following basic fact: 
 
 \begin{lem} \label{nitpick}
 If $M$ is a a smooth compact oriented $4$-manifold with $b_+ > 1$,
 and if $g$ is a K\"ahler-Einstein metric on $M$ with negative 
 scalar curvature, then equality is achieved in \eqref{mocha} by $g$. 
 \end{lem}
 \begin{proof}
For any compact K\"ahler surface $(M,J)$ with $b_+>1$, the 
classical Seiberg-Witten invariant is well-defined and non-zero \cite{witten}  for the spin$^c$
structure determined by $J$, and ${\mathbf a}= c_1^\RR (M, J)$
is therefore a monopole class. Hence 
$c_1^\RR (M, J)\in {\mathfrak C} \subset \hull ({\mathfrak C})$, and
$$
\beta^2 (M) = \max \{ {\mathbf v}^2 ~|~ {\mathbf v} \in \hull ({\mathfrak C})\} 
\geq c_1^2 (M) . 
$$
On the other hand, the Ricci form $\rho = r (J\cdot , \cdot )$ represents
$2\pi c_1^\RR (M, J)$ in de Rham cohomology, and just equals 
 $s\omega/4$
in the  K\"ahler-Einstein case. Thus, since the 
volume form of a K\"ahler surface is given by $\omega^2/2$, we have 
$$
\int_M s^2d\mu =  \int  \frac{(s\omega)^2}{2}\\
                             = 8 \int \rho \wedge \rho \\
                              =  32\pi^2 c_1^2 (M) ~.
$$
 Proposition \ref{best} therefore tells us that  
$$ 32\pi^2  c_1^2 (M) = \int_M s^2 d\mu \geq 32\pi^2 \beta^2 (M) \geq 32\pi^2  c_1^2 (M), $$
and equality must thus hold at every step. Hence $\beta^2 (M)= c_1^2(M)$, and
equality is achieved in \eqref{mocha} by $g$, as claimed. 
 \end{proof}

 \begin{lem} \label{nitwit}
 Let $M$ be a compact oriented $4$-manifold with $b_+\geq 2$
which  carries a non-zero  monopole class. Then 
whenever equality holds in  \eqref{mocha}
for a 
metric 
$g$ on $M$, equality holds
in \eqref{java}, too. 
 \end{lem} 
 \begin{proof}
 If equality holds in  \eqref{mocha}, we have
 $$
32\pi^2 \beta^2 (M) = \int_M s^2 d\mu ~, 
$$
so  by Propositions \ref{best} and  \ref{biga}, there is a 
monopole class ${\mathbf a}$ such that 
$$
32\pi^2 \beta^2 (M) = \int_M s^2d\mu \geq
{32\pi^2}({\mathbf a}^+)^2\geq {32\pi^2}\beta^2 (M) > 0, 
$$
and  equality must therefore hold throughout. 
But Proposition \ref{best} then asserts that there exists 
a complex structure
$J$  such that $(M,g,J)$ is a 
a K\"ahler manifold of constant negative scalar curvature. 

Now  
any K\"ahler metric on a complex surface automatically satisfies 
$|W_+|^2 = s^2/24$, so that $s-\sqrt{6}|W_+|= \frac{3}{2} s$ wherever $s\leq 0$.
Our negative-scalar-curvature K\"ahler metric $g$ thus satisfies 
$$
\int_M (s-\sqrt{6}|W_+|)^2d\mu = \left(\frac{3}{2}\right)^2 \int_M s^2d\mu = 72\pi^2 \beta^2 (M) ,
$$
and  therefore also achieves equality  in \eqref{java}, as claimed.
 \end{proof}

We  now analyze  the boundary case of \eqref{java}. 

\begin{thm} \label{newnew} 
Let $M$ be a compact oriented $4$-manifold with $b_+\geq 2$
which  carries a non-zero  monopole class, and suppose that 
$g$ is a  metric  on $M$ such that equality holds in \eqref{java}: 
$$
\int_M (s-\sqrt{6}|W_+|)^2d\mu = 72\pi^2 \beta^2(M).
$$
Then $g$ is K\"ahler-Einstein, with negative Einstein constant. 
\end{thm}
\begin{proof}
Let ${\mathbf v}\in \hull({\mathfrak C})$ be a point where ${\mathbf v}^2={\mathbf v}\cdot {\mathbf v}$
achieves its maximum value, namely $\beta^2(M)$.
Let ${\mathbf a}_1, \ldots , {\mathbf a}_n\in {\mathfrak C}$
be a list of all the monopole classes, and express ${\mathbf v}\in \hull ({\mathfrak C})$
as 
$${\mathbf v} = \sum_{j=1}^n t_j {\mathbf a}_j$$
where the coefficients $t_j \in [0,1]$ satisfy $\sum_j t_j =1$; and 
after permuting the ${\mathbf a}_j$ as necessary,  we may henceforth assume that   
$t_j> 0$  iff $j\leq m$, where $m$ is some integer, $1\leq m \leq n$. 
By Propositions \ref{upscale} and  \ref{biga},
\begin{eqnarray*}
\frac{1}{72\pi^2}\int_M (s-\sqrt{6}|W_+|)^2d\mu &\geq&
 \max ~\{  ({\mathbf a}^+_j)^2~|~ j=1, \ldots , n\} \\
 &\geq&  ({\mathbf v}^+)^2\geq {\mathbf v}^2= \beta^2 (M) 
\end{eqnarray*}
and  our hypotheses therefore imply that equality holds at every step. In particular, it
follows that ${\mathbf v}={\mathbf v}^+$ and that $ \max_j ({\mathbf a}^+_j)^2= \beta^2 (M)$. 
Since the intersection form is positive definite
on  ${\mathcal H}^+_g$, the Cauchy-Schwarz inequality therefore tells us that 
$$
{\mathbf v}\cdot {\mathbf a}_j^+ \leq \sqrt{ ({\mathbf a}^+_j)^2} \sqrt{ ({\mathbf v})^2}\leq \beta^2(M)~,
$$
for all $j$, with equality iff ${\mathbf a}_j^+= {\mathbf v}$. Since 
\begin{eqnarray*}
\beta^2 (M) &=&{\mathbf v}\cdot {\mathbf v}={\mathbf v}\cdot {\mathbf v}^+
\\&=& {\mathbf v}\cdot  \left(\sum_{j=1}^m t_j {\mathbf a}_j^+\right) \\
&=&   \sum_{j=1}^m t_j \left({\mathbf v}\cdot {\mathbf a}_j^+\right)\\
&\leq& \sum_{j=1}^m t_j ~\beta^2 (M)\\ 
&=& \beta^2 (M) \left(\sum_{j=1}^m t_j \right) \\
&=& \beta^2 (M) ~, 
\end{eqnarray*}
we must therefore have 
${\mathbf a}_j^+ ={\mathbf v}$   for every $j= 1, \ldots , m$. 

For each $j= 1, \ldots , m$, we therefore have $({\mathbf a}_j^+)^2=\beta^2(M)$.
Moreover, $\beta^2(M) > 0$ by Lemma \ref{betapos}. Our hypotheses  thus imply that  
$$\int (s-\sqrt{6}|W_+|)^2d\mu = 72\pi^2 ({\mathbf a}_j^+)^2 > 0,$$
and Proposition \ref{upscale} therefore tells us that 
 there  is a $g$-compatible symplectic form 
$\omega_j$ such that  $[\omega_j]$ is a negative multiple of ${\mathbf a}_j^+= {\mathbf v}$,
and such that  
$c_1^\RR (M,\omega_j) = {\mathbf a}_j$ for each $j=1, \ldots , m$. 
Since $[\omega_j]^2/2= \text{Vol} (M,g)$ for each $j$, it follows that 
$[\omega_1] = \cdots = [ \omega_m]\in H^2 (M, \RR)$. 
But  each  $\omega_j$ is harmonic with respect to $g$, 
and the harmonic representative of any de Rham class is unique. 
Hence 
$\omega_1= \cdots = \omega_m$. But since $c_1^\RR (M,\omega_j) = {\mathbf a}_j$, 
this implies  that $ {\mathbf a}_1= \cdots =  {\mathbf a}_m$.
Hence 
$m=1$, and
$${\mathbf v}= \sum_{j=1}^m t_j {\mathbf a}_j = {\mathbf a}_1 = c_1^\RR (M, \omega ).$$ 

Let us now simplify our notation by setting $\omega= \omega_1$. 
Since  $- [\omega]\propto  {\mathbf v}= c_1(M, \omega )$, 
the curvature  of any connection on 
the anti-canonical line bundle $L$ of $(M, \omega)$ 
must be cohomologous to a  constant negative multiple  of  $\omega$.  However, 
 we saw in (\ref{bla1}--\ref{bla2}) that 
the curvature $F_B=F_B^++F_B^-$ of the Blair connection on $L$ 
is given by 
\begin{eqnarray*}
iF_B^+ &=& \frac{s+s^*}{8}\omega+  W^+(\omega )^\perp   \\
iF_B^-&=&  \frac{s-s^\star }{8}\hat{\omega}+ \mathring{\varrho}\end{eqnarray*}
where 
  $W^+(\omega )^\perp$
  is the component of $W^+(\omega )$ orthogonal to $\omega$,  
  $$\mathring{\varrho}(\cdot , J\cdot )  = \frac{\mathring{r} +  J^*\mathring{r} }{2} ,$$
and where  the bounded anti-self-dual $2$-form 
 $\hat{\omega}\in \Lambda^-$
 is  defined only on the open set where $s^\star - s\neq 0$,  and  satisfies 
 $|\hat{\omega}|\equiv \sqrt{2}.$
 Here, the star-scalar curvature $s^*$  once again means the important quantity 
 $$s^{*}=s+ |\nabla \omega |^2 = 2W_+ (\omega , \omega ) +  \frac{s}{3}. $$
Since  Proposition \ref{upscale} tells us that 
$(M,g,\omega )$ is saturated, $s+s^*$ is constant and $W^+(\omega )^\perp= 0$. 
Hence  $F_B^+$  is closed, and therefore   $\star F_B=2F_B^+-F_B$ 
is closed, too. Thus $F_B$ is harmonic.  
But we  also know that $F_B$ is cohomologous to a constant multiple of $\omega$,
which 
 is itself a self-dual harmonic form. Hence 
 $F_B^-\equiv 0$, and 
$$\mathring{\varrho}\equiv  \frac{s^\star -s}{8}\hat{\omega}.$$
This shows that  
$$|\mathring{r}|^2 \geq \frac{(s^*-s)^2}{16}$$
at every point of $M$, with equality precisely at those points at which 
the Ricci tensor $r$ is $J$-invariant. 

On the other hand, $W_+$ has eigenvalues $(-\lambda/2, -\lambda/2, \lambda)$,
where
$$\lambda = \frac{1}{2} W_+(\omega , \omega )= \frac{3s^*-s}{12} ,$$
so 
$$|W_+|^2 = \frac{(3s^*-s)^2}{96}.$$
Hence
\begin{eqnarray*}
4\pi^2 (2\chi + 3\tau ) (M)  & = & \int_M \left(\frac{s^2}{24}+ 2 |W_+|^2 - \frac{|\mathring{r}|^2}{2}\right)d\mu  \\
 & = &  \int_M   \left(\frac{s^2}{24}+  \frac{2(3s^*-s)^2}{96}- \frac{|\mathring{r}|^2}{2}\right)d\mu \\
  & \leq &    \int_M \left(\frac{s^2}{24}+  \frac{2(3s^*-s)^2}{96}- \frac{(s^*-s)^2}{32}\right)d\mu \\
   & = &  \frac{1}{32}  \int_M  \left(s^2-2ss^* +5(s^*)^2\right)d\mu 
\end{eqnarray*}
with equality iff $|\mathring{r}|^2 \equiv (s^*-s)^2/16$.
On the other hand, since $F_B=F_B^+$,  
\begin{eqnarray*}
4\pi^2 (2\chi + 3\tau ) (M)  & = & 4\pi^2 c_1^2 (M)   \\
 & = &  \int_M  (\frac{s+s^*}{8}\omega )\wedge  (\frac{s+s^*}{8}\omega ) \\
  & =&  \frac{1}{32}   \int_M  \left(s^2+2ss^* +(s^*)^2\right)d\mu 
\end{eqnarray*}
so we therefore have
$$
\int_M  \left(s^2-2ss^* +5(s^*)^2\right)d\mu \geq   \left(s^2+2ss^* +(s^*)^2\right)d\mu , 
$$
which we can rewrite as 
\begin{equation}
\label{gizmo}
\int_M  4s^*(s^*-s)  d\mu \geq  0 ~; 
\end{equation}
moreover,   equality can hold only if  $|\mathring{r}|^2 \equiv (s^*-s)^2/16$.
However, since $(M,g,\omega)$ is saturated,  $s^*+s$  is a negative constant,
and $W_+(\omega , \omega )  \leq 0$; hence 
$s^*\leq s/3$, and  $s^* \leq (s+s^*)/4$ is therefore negative everywhere. 
 On the other hand, 
 $s^*-s = |\nabla \omega |^2 \geq 0$  on  any almost-K\"ahler  manifold.
 Hence 
$$s^* (s^*-s) \leq 0$$
everywhere on $M$, 
with equality only at points where $s=s^*$. 
The  inequality (\ref{gizmo}) therefore  implies that 
 $$|\nabla \omega |^2 = s-s^*\equiv 0.$$ Hence  $(M,g,\omega)$ is K\"ahler.
But  equality in (\ref{gizmo}) only holds  if $|\mathring{r}|^2 \equiv (s^*-s)^2/16$, so 
we moreover must have $\mathring{r} \equiv 0$, and 
we therefore conclude  that $(M,g)$ is K\"ahler-Einstein, as promised. 
\end{proof}

Our main result now  follows easily:

\begin{proofA} Theorem \ref{magna}
  shows that \eqref{mocha} and \eqref{java} hold for any metric on any 
  $4$-manifold with $b_+\geq 2$. On the other hand, assuming there is at
  least one non-zero monopole class, 
  Theorem \ref{newnew} shows that any  metric for which 
  equality holds in \eqref{java} must be K\"ahler-Einstein. 
  Lemma \ref{nitwit} thus implies that any  metric for which 
  equality holds in \eqref{mocha} must be K\"ahler-Einstein, too. 
  Finally, Lemmas \ref{nitpick} and \ref{nitwit} show that 
  equality actually does  hold  in \eqref{mocha} and \eqref{java} 
  when the metric is 
  K\"ahler-Einstein. 
  \end{proofA}

Of course, the method  used here to treat  the 
boundary case of  \eqref{mocha} 
has a   Rube Goldberg feel to it, since  it
proceeds by  reducing  an easy 
problem to  a harder one. However, 
it is not difficult to  
 winnow 
a simple,  direct treatment of this case 
 out of  the  above 
discussion.   Details are left to  
  the interested reader.

\section{Concluding Remarks} 

One apparent weakness of our definition of $\beta^2(M)$ is that 
there is no  obvious way of exactly  determining the entire set ${\mathfrak C}(M)$ of all 
 monopole classes of a given $4$-manifold $M$. However, we {\em do} have 
 various criteria which serve to show that certain classes really do belong to 
 ${\mathfrak C}(M)$.  Thus, if ${\mathfrak S}\subset {\mathfrak C}(M)$ is some 
 collection of known monopole classes,  we then have 
$$
\beta^2 (M) \geq \max \{ {\mathbf v}^2~|~ {\mathbf v} \in \hull ( {\mathfrak S}) \} . 
$$
It is thus relatively easy 
to find lower bounds for $\beta^2$, even without knowing ${\mathfrak C}(M)$
exactly. 

At the same time, our curvature estimates \eqref{mocha} and \eqref{java} 
provide   upper bounds for $\beta^2 (M)$ for each metric 
$g$ on $M$. By taking an infimum of such upper bounds for a
carefully chosen  sequence of metrics $g_j$ on $M$, one can, 
in practice, often  
 determine $\beta^2(M)$ by  showing that it is simultaneously 
no less than and no greater than  some target value. 

\begin{xpl}
Let $X$ be a minimal complex surface of general type, and 
let $M= X\# k \overline{\bcp}_2$ be its blow-up at $k$ points. 
Then 
$$\beta^2 (M)  = c_1^2 (X).$$
Indeed, if $E_1, \ldots , E_k$ are generators for the 
various copies of $H^2 (\overline{\bcp}_2, \ZZ)$,  
then $\pm c_1 (X) \pm E \pm \cdots \pm E_k$
are the first Chern classes of various complex structures of K\"ahler type
on $M$, and so are monopole classes \cite{witten}. 
Hence $c_1(X) \in \hull ({\mathfrak C}(M))$, and hence 
$\beta^2 (M) \geq c_1^2 (X)$. However, by approximating
the K\"ahler-Einstein orbifold metric on the pluricanonical model for
$X$, one can construct \cite{lno} 
sequences
of metrics $g_j$ on $M$ with $\int s^2 d\mu \searrow 32\pi^2 c_1^2 (X)$.
Thus 
\eqref{mocha} implies  that 
we 
also have $c_1^2 (X) \geq \beta^2 (M)$, and the claim follows. 
\end{xpl}

\begin{xpl} \label{hiho} 
Let $X$, $Y$, and $Z$ be simply connected, minimal complex surfaces of general type
with $h^{2,0}$ odd. Let $M= X\# Y \# Z\# k \overline{\bcp}_2$.
Then 
$$
\beta^2 (M ) = c_1^2(X)+  c_1^2(Y) + c_1^2 (Z).
$$
Indeed, using the Bauer-Furuta invariant, one can show that 
$$\pm c_1(X) \pm c_1 (Y) \pm c_1(Z)  \pm E_1  \pm \cdots  \pm E_k
\in {\mathfrak C}(M).$$ Hence ${\mathbf v} = c_1(X)+  c_1(Y) + c_1 (Z) \in
\hull ({\mathfrak C}(M))$, and 
$$\beta^2 (M ) \geq [c_1(X)+  c_1(Y) + c_1 (Z)]^2=   c_1^2(X)+  c_1^2(Y) + c_1 (Z).$$ 
On the other hand, there exist \cite{il2} sequences
of metrics $g_j$ on $M$ with $\int s^2 d\mu \searrow 32\pi^2 [c_1^2(X)+  c_1^2(Y) + c_1 (Z)]$, so 
\eqref{mocha} therefore
shows that 
we 
also have $c_1^2(X)+  c_1^2(Y) + c_1^2 (Z)\geq \beta^2 (M)$. The claim therefore follows. 

Similar techniques can also be used for connected sums involving  two or four 
surfaces of general type. 
\end{xpl}

\begin{xpl}
Let $N$ be any oriented $3$-manifold, and let $M= N\times S^1$. Then 
$\beta^2 (M) =0$, because one has  $\int s^2 d\mu \searrow 0$ for product
metrics on $M$ with shorter and shorter $S^1$ factors. However, that 
results of Kronheimer \cite{K} imply  that such manifolds  typically carry
many monopole classes,  although these  all belong to the 
isotropic subspace $H^2(N) \hookrightarrow H^2 (N\times S^1)$.
\end{xpl}

By the arguments detailed in \cite{lmo,lric}, the estimates 
\eqref{mocha} and \eqref{java} have the following interesting consequences: 

\begin{thm}
Let $M$ be a smooth compact oriented $4$-manifold with $b_+(M) \geq 2$.
If $M$ admits an  Einstein metric $g$, then 
$$
(2\chi - 3\tau ) (M) \geq \frac{1}{3} \beta^2 (M).
$$
Moreover, if $M$ carries a non-zero monopole class, 
equality occurs 
only if $(M,g)$ is a compact quotient ${\mathbb C}{\mathcal H}_2/\Gamma$
of the complex hyperbolic plane, equipped with 
a constant multiple of its standard K\"ahler-Einstein metric.
\end{thm}

\begin{thm}
Let $M$ be a smooth compact oriented $4$-manifold with $b_+(M) \geq 2$.
If $M$ admits an  Einstein metric $g$, then 
$$
(2\chi + 3\tau ) (M) \geq \frac{2}{3} \beta^2 (M),
$$
with equality only if both sides vanish, in which case 
$g$
must be  a hyper-K\"ahler metric, and 
$M$ must be  diffeomorphic to $K3$ or $T^4$. 
\end{thm}

\begin{thm}
Let $M$ be a smooth compact oriented $4$-manifold with $b_+(M) \geq 2$.
Then any metric $g$ on $M$  satisfies
$$
\int_M |r|^2 d\mu \geq 8\pi^2 \left[ 2\beta^2 - (2\chi + 3\tau )\right] (M) , 
$$
with equality iff $g$ is K\"ahler-Einstein. 
\end{thm}

Now Proposition \ref{yobo} entitles us to introduce the following definition: 

\begin{defn}
If $M$ is any  smooth compact oriented $4$-manifold with $b_+(M) \geq 2$,
we set $\beta (M) := \sqrt{\beta^2 (M)}$.  
\end{defn}

This invariant provides a natural yardstick with which to measure  the Yamabe invariants
of  $4$-manifolds:

\begin{thm} \label{yammer} 
Let $M$ be a smooth compact oriented $4$-manifold with $b_+(M) \geq 2$.
If $M$ carries at least one non-zero monopole class, then 
the Yamabe invariant of $M$ satisfies 
$$
{\mathcal Y} (M) \leq - 4\sqrt{2} \pi  \beta (M) .
$$
\end{thm}

We remark in passing that if $M$ does not admit a metric of positive scalar curvature,
its Yamabe invariant ${\mathcal Y}(M)$ is just the supremum of the 
scalar curvatures of unit-volume constant-scalar-curvature metrics on 
$M$. This result is therefore an immediate consequence of  \eqref{mocha}. 
More intriguingly, though,  Theorem \ref{yammer}  is actually sharp; equality actually  holds 
 \cite{lno,il2} for large classes of 
$4$-manifolds, including those discussed on page \pageref{hiho}.

Now, while we have seen that considering the convex hull
of the set of monopole classes leads to an elegant invariant $\beta^2(M)$ 
which seems remarkably well adapted to the study of the curvature
of $4$-manifolds, it is still unclear whether this approach is optimal 
in all circumstances. Indeed, the basic forms of our estimates,
seen in Propositions \ref{best} and \ref{upscale}, involve the
numbers $({\mathbf a}^+)^2$ for the various monopole classes, and 
one can therefore \cite{lebsurv2} define an invariant which 
simply tries to  make optimal use of this information. 
Indeed, 
consider the open Grassmannian ${\mathbf G \mathbf r} = Gr^+_{b_+}[H^2(M, \RR )]$ of  all 
maximal linear subspaces ${\mathbf H}$  of the 
second cohomology on which the restriction of the intersection pairing
  is positive definite. Each 
element ${\mathbf H}\in {\mathbf G \mathbf r}$ then determines 
an orthogonal  decomposition 
$$H^2(M, \RR ) = {\mathbf H} \oplus {\mathbf H}^\perp$$
with respect to the intersection form.  Given a monopole class ${\mathbf a}\in {\mathfrak C}$
and a positive subspace ${\mathbf H}\in {\mathbf G \mathbf r}$, we may then  define 
${\mathbf a}^+$ to be  the orthogonal projection  of ${\mathbf a}$ 
into $\mathbf H$. Using this, we now  define
yet another  oriented-diffeomorphism invariant.
\begin{defn}
Let $M$ be a smooth compact oriented manifold with $b_+\geq 2$.
If ${\mathfrak C}=\varnothing$, set  $\alpha^2 (M)=0$.
Otherwise, we set
$$
\alpha^2 (M) = \inf_{{\mathbf H}\in {\mathbf G \mathbf r}}\left[ \max_{a \in 
{\mathfrak C}} ~{({\mathbf a}^+)^2} \right] ~ .
$$
\end{defn}
Propositions \eqref{best} and \eqref{upscale} then 
easily imply that \eqref{mocha} and \eqref{java} 
still hold when $\beta^2(M)$ is replaced by $\alpha^2 (M)$. Moreover,
the proof of Proposition \ref{biga} shows that one always has 
$$\alpha^2 (M) \geq \beta^2(M).$$

On the other hand, we have also seen that \eqref{mocha} and \eqref{java} 
are sharp for large classes of manifolds, such as those  
discussed on page \pageref{hiho}. Thus $\alpha^2 = \beta^2$ in all
these cases. It is therefore only natural for 
 us to ask whether this is a general phenomenon. 
In this direction, however, we can only give some  partial results. 
We begin with  the following:

\begin{lem}\label{beggar} 
Let $M$ be a smooth oriented $4$-manifold with $b_+\geq 2$.
Then 
$$\alpha^2 (M) =0 \Longleftrightarrow \beta^2 (M) =0.$$
\end{lem}
\begin{proof} The $\Longrightarrow$  direction is obvious, since $\alpha^2 \geq \beta^2 \geq 0$.
Conversely, if $\beta^2 =0$, the intersection form must be negative-semi-definite
on $\span ({\mathfrak C})$. Write this subspace as ${\mathbf N} \oplus {\mathbf I}$,
where the intersection form is negative-definite on ${\mathbf N}$ and vanishes
on ${\mathbf I}$. We can then choose a sequence ${\mathbf H}_j \in {\mathbf G \mathbf r}$
which are all orthogonal to ${\mathbf N}$ and which decompose orthogonally as  
${\mathbf P} \oplus {\mathbf J}_j$, where ${\mathbf P}$ is orthogonal to ${\mathbf I}$
and ${\mathbf J}_j\to {\mathbf I}$. Then each monopole class satisfies $({\mathbf a}^+)^2\to 0$
for  this sequence. It thus follows that $\alpha^2=0$, as claimed.  
\end{proof}

Next, we point out the following: 

\begin{prop} \label{couscous} 
Let $M$ be a smooth oriented $4$-manifold with $b_+\geq 2$. Suppose, moreover, that there 
is a linear subspace 
${\mathbf L} \subset H^2 (M, \RR)$
on which the intersection form is of  Lorentzian $({+}{-}\cdots {-})$  type, with 
$${\mathfrak C}(M) \subset {\mathbf L} \subset H^2 (M, \RR).$$
Then $\alpha^2(M) =\beta^2 (M)$. 
\end{prop}
\begin{proof}By Lemma \ref{beggar}, we may assume that $\beta^2 (M )> 0$. Let
${\mathbf v}\in \hull ({\mathfrak C}) \subset {\mathbf L}$ be an 
element with  ${\mathbf v}^2=\beta^2 (M) > 0$. Now, since 
$(1-t) {\mathbf v} + t {\mathbf a}\in \hull ({\mathfrak C})$ for any ${\mathbf a }\in {\mathfrak C}$ and any
$t\in [0,1]$, we therefore have 
$${\mathbf v}^2\geq [(1-t) {\mathbf v} + t {\mathbf a}]^2 = {\mathbf v} ^2 + 2t (
 {\mathbf v}\cdot {\mathbf a}- 
{\mathbf v}^2) + O(t^2)$$
for all small positive $t$, and it therefore follows that 
$$ {\mathbf v}\cdot {\mathbf a} \leq {\mathbf v}^2$$
for all monopole classes ${\mathbf a}$. Since ${\mathfrak C}(M)$ is 
invariant under multiplication by $-1$, it moreover follows that 
$$
| {\mathbf v}\cdot {\mathbf a} |\leq {\mathbf v}^2 ~~\forall {\mathbf a}\in {\mathfrak C}(M).
$$

 Now let ${\mathbf P}\subset {\mathbf L}^\perp$
be a maximal positive subspace, and set ${\mathbf H}= {\mathbf P}\oplus \span ({\mathbf v})$.
Then for this choice of  ${\mathbf H}\in {\mathbf G \mathbf r}$ we have 
$${\mathbf a}^+ = \frac{{\mathbf v}\cdot {\mathbf a}}{ {\mathbf v}^2} {\mathbf v}$$
and hence 
$$({\mathbf a}^+)^2 =  
\frac{({\mathbf v}\cdot {\mathbf a})^2}{{\mathbf v}^2}\leq {\mathbf v}^2=
\beta^2(M) $$
for all ${\mathbf a}\in {\mathfrak C}$. Hence
$$\alpha^2 (M) = \inf_{{\mathbf H}\in {\mathbf G \mathbf r}}\left[ \max_{a \in 
{\mathfrak C}} ~{({\mathbf a}^+)^2} \right] \leq \beta^2 (M).$$
But we also know that $\beta^2\leq \alpha^2$, so it follows that 
  $\alpha^2=\beta^2$, as claimed. 
\end{proof}

\begin{xpl} If $(M,J)$ is a compact complex surface of 
K\"ahler type with $b_+ > 1$, we may take ${\mathbf L} = H^{1,1}(M, \RR)\subset H^2 (M, \RR)$.
Since an argument due to Witten \cite{witten} shows that  solutions of the Seiberg-Witten 
equations  can  exist with respect to a K\"ahler metric only when 
$c_1(L)$ is a $(1,1)$-class,  it follows that 
any monopole class must belong to  ${\mathbf L}$.
This provides one explanation of why  $\alpha^2=\beta^2$ for complex algebraic surfaces. 
\end{xpl}

In light of Proposition \ref{couscous}, the reader may be curious as to why we 
have systematically excluded the case of $b_+=1$ in this paper.
In truth, most of the formal theory actually works perfectly well in this case. 
However, the Seiberg-Witten invariants have a chamber structure when
$b_+=1$, and this has the effect that, for example, complex surfaces
with $c_1^2< 0$ will typically not carry any monopole classes at all.
Nonetheless,  Seiberg-Witten theory still gives rise \cite{lno} to non-trivial
curvature bounds in this setting, even though this phenomenon  cannot be 
explained in terms of    monopole classes. 

We now turn to a more complicated situation: 

\begin{prop}
Let $M$ be a smooth oriented $4$-manifold with $b_+\geq 2$, and suppose that there 
is a collection of mutually orthogonal linear subspaces 
${\mathbf L}_j \subset H^2 (M, \RR), j=1, \ldots , \ell,$
on each of which the intersection form is of  Lorentzian $({+}{-}\cdots {-})$  type. 
Moreover, suppose that  
$${\mathfrak C} (M)= {\mathfrak C}_1 \times \cdots \times {\mathfrak C}_\ell
\subset {\mathbf L}_1 \oplus \cdots \oplus {\mathbf L}_\ell, $$
for some subsets
$${\mathfrak C}_j \subset {\mathbf L}_j, ~~~  j=1, \ldots , \ell .$$
Then $\alpha^2(M) =\beta^2 (M)$. 
\end{prop}
\begin{proof}
Fix a maximal positive subspace ${\mathbf P}\subset ({\mathbf L}_1 
\oplus \cdots \oplus {\mathbf L}_\ell)^\perp$, and consider choices of 
${\mathbf H}\in {\mathbf G \mathbf r}$ of the form 
${\mathbf H}= {\mathbf P}\oplus \span \{ e_1 , \ldots , e_\ell\}$,
where $e_j \in {\mathbf L}_j$ is a non-zero time-like  vector. 
If the intersection form on $\span ({\mathfrak C}_j)$ is negative-definite, 
moreover choose  $e_j\in {\mathbf L}_j$ to  be orthogonal to this subspace. If, on the other hand,  
$\span ({\mathfrak C}_j)$ is 
Lorentzian, set $e_j={\mathbf v}_j$, where ${\mathbf v}_j$ maximizes
${\mathbf v}^2$ on $\hull ({\mathfrak C}_j)$. Finally, if the intersection form is 
degenerate on $\span ({\mathfrak C}_j)$, choose  
${\mathbf v}_j\in \hull ({\mathfrak C}_j)$ to be a non-zero null vector, 
and consider a sequence of different possible 
$e_j$ converging to ${\mathbf v}_j$.  In this way, one 
obtains a sequence of choices of  ${\mathbf H}$ for which 
$\max ({\mathbf a}^+)^2 \to \sum ({\mathbf v}_j)^2 = \beta^2 (M)$.
Hence $\alpha^2 \leq \beta^2 \leq \alpha^2 $, and 
$\alpha^2 (M) = \beta^2 (M)$, as claimed. 
\end{proof}

This result gives a partial explanation of why $\alpha^2 = \beta^2$
for the connected sums of complex surfaces we have considered, 
since  the set of  known monopole classes in this case 
constitutes a configuration of the described type, where the 
Lorentzian subspaces in question are given by 
$H^{1,1}$ of the various summands. Of course, this
explanation still remains less than entirely satisfactory, since 
we cannot be  absolutely certain that we currently have a full catalog  of the monopole 
classes of these spaces. 

Finally, let us point out that  one {\em cannot} hope to prove 
that $\alpha^2=\beta^2$  if ${\mathfrak C}$ is simply replaced with an 
arbitrary finite, centrally symmetric set in an arbitrary finite-dimensional 
vector space with indefinite inner product. For example, let us just 
consider $\RR^3$ equipped with the $({+}{+}{-})$ inner product 
$dx^2 + dy^2 -dz^2$, and consider the candidate for 
``${\mathfrak C}$'' given by 
$$\left\{ \pm (1,0,1), ~\pm (\frac{\sqrt{3}}{2}, -\frac{1}{2} , 1),  ~\pm  (-\frac{\sqrt{3}}{2}, -\frac{1}{2} , 1)
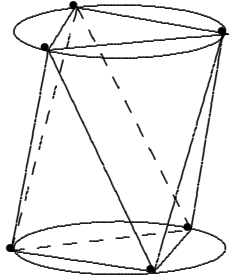
\right\}, 
$$
\begin{center}
\mbox{
\beginpicture
\setplotarea x from 40 to 150, y from  -20 to 100
\ellipticalarc axes ratio 4:1  360 degrees from  140 85
center at   100 85
\ellipticalarc axes ratio 4:1  360 degrees from  60 3
center at   100 3
\put {\circle*{3}} [B1] at 140  85
\put {\circle*{3}} [B1] at  73  79
\put {\circle*{3}} [B1] at  84  95
\put {\circle*{3}} [B1] at  60  3
\put {\circle*{3}} [B1] at  113  -5
\put {\circle*{3}} [B1] at  127  11
{\setlinear 
\plot  141 84  73 78 /
\plot  141 84  84 94 /
\plot  83 94  73  78 /
\plot   59 2 113  -6  /
\plot 127 10  113 -6  /
\plot   59 2  73  78  /
\plot  127 10 139 84  /
\plot   113  -6  139 84 /
\plot   113  -6  73 78 /
\setdashes 
\plot  127 10   84 94 /
\plot   59  2  127 10 /
\plot   59  2  84 94 /
}
\endpicture
}
\end{center}
Because the elements of this configuration ``${\mathfrak C}$'' are all null vectors, one can use 
Proposition \ref{couscous} 
``upside-down'' to 
show that ``$\alpha^2$'' must then equal $1$. On the other hand, 
a simple symmetry argument shows that ``$\beta^2$'' equals $\frac{3}{4}$ for
this configuration. Of course, this this choice of  ``${\mathfrak C}$''  does not consist
of integer points, but one can easily  remedy this by  rational approximation and
rescaling.  

The upshot is that while one definitely has $\alpha^2 (M)=\beta^2(M)$ for a 
remarkably large and interesting array of examples, this statement 
can generally  hold  only to the  degree that 
the set ${\mathfrak C}$ of monopole classes 
satisfies some manifestly 
non-trivial geometric 
constraints.
The precise extent to  which these constraints hold or fail 
remains to be determined. It is  hoped that  
some  interested reader will find the  challenge 
of  fully  fathoming this mystery
both stimulating and fruitful.

  \end{document}